\documentclass[12pt,reqno]{amsart}
\usepackage{amsthm}
\usepackage{cases}
\usepackage{amscd}
\usepackage{amsfonts}

\usepackage{amssymb}
\usepackage{amsmath}
\usepackage[dvipdfmx]{graphicx}
\usepackage{booktabs}
\usepackage{adjustbox}
\usepackage{subcaption}
\usepackage[percent]{overpic}
\usepackage{hyperref}
\usepackage{graphicx}
\usepackage{bm}
\usepackage[utf8]{inputenc}   
\usepackage[T1]{fontenc}      
\usepackage{textcomp}         
\usepackage{newunicodechar}
\newunicodechar{≥}{\ensuremath{\geq}}
\newunicodechar{≤}{\ensuremath{\leq}}
\newunicodechar{–}{\textendash}      
\newunicodechar{：}{:}
\DeclareUnicodeCharacter{2212}{\ensuremath{-}}

\usepackage{algorithm}
\usepackage{algpseudocode} 
\renewcommand{\algorithmicrequire}

\usepackage{cleveref}
\usepackage{color}
\usepackage{datetime}

\allowdisplaybreaks

\makeatletter
\newcommand\figcaption{\def\@captype{figure}\caption}
\newcommand\tabcaption{\def\@captype{table}\caption}
\makeatother

\usepackage{geometry}
\usepackage{algorithm}
\usepackage{multirow}
\usepackage{mathrsfs}
\usepackage{fancyhdr}
\newcommand{\beq}{\begin{equation}}
\newcommand{\eeq}{\end{equation}}
\newcommand{\bea}{\begin{eqnarray}}
\newcommand{\eea}{\end{eqnarray}}
\newcommand{\beas}{\begin{eqnarray*}}
\newcommand{\eeas}{\end{eqnarray*}}

\begin{document}
\title{Robust PDE Discovery under Sparse and Highly Noisy Conditions via Attention Neural Networks}

\author[S.L. Zhang, Y.Q. Huang, N.Y. Yi , S.H. Zhang]{Shilin Zhang$^\dagger$, Yunqing Huang$^{\ddagger,*}$, Nianyu Yi$^\S$ and Shihan Zhang$^\dagger$}

\address{$^\dagger$ School of Mathematics and Computational Science, Xiangtan University, Xiangtan 411105, P.R.China }\email{shilinzhang@smail.xtu.edu.cn (S.L. Zhang), shihanzhang@smail.xtu.edu.cn (S.H. Zhang)}
\address{$\ddagger$ National Center for Applied Mathematics in Hunan, Key Laboratory of Intelligent Computing \& Information Processing of Ministry of Education, Xiangtan University, Xiangtan 411105, Hunan, P.R.China} \email{huangyq@xtu.edu.cn}
\address{$^\S$ Hunan Key Laboratory for Computation and Simulation in Science and Engineering; School of Mathematics and Computational Science, Xiangtan University, Xiangtan 411105, P.R.China}
\email{yinianyu@xtu.edu.cn}



\begin{abstract}

The discovery of partial differential equations (PDEs) from experimental data holds great promise for uncovering predictive models of complex physical systems. In this study, we introduce an efficient automatic model discovery framework, ANN-PYSR, which integrates attention neural networks with the state-of-the-art PySR symbolic regression library. Our approach successfully identifies the governing PDE in six benchmark examples. Compared to the DLGA framework, numerical experiments demonstrate ANN-PYSR can extract the underlying dynamic model more efficiently and robustly from sparse, highly noisy data (noise level up to 200\%, 5000 sampling points). It indicates an extensive variety of practical applications of ANN-PYSR, particularly in conditions with sparse sensor networks and high noise levels, where traditional methods frequently fail.
\end{abstract}


\maketitle


\section{Introduction}

Through long-term astronomical observations and geometric intuition, Johannes Kepler formulated his third law of planetary motion, a fundamental principle that not only unraveled celestial dynamics but also inspired Isaac Newton's universal gravitation theory. Similarly, Planck’s law of radiation, derived by fitting data, laid the foundation for quantum mechanics. These historic examples underscore the power of integrating data with theory to advance scientific discovery. For centuries, scientists have sought to model physical phenomena using partial differential equations (PDEs) derived from first principles, which form the foundation of classical mechanics and relativity \cite{1}. However, the modeling of contemporary complex systems, such as power grids, climate dynamics, and biological processes, has exposed the limitations of purely first-principles approaches. In climate modeling, for example, the challenge of unobservable variables (e.g. deep ocean temperatures) and errors introduced by experimental noise in derivative calculations can severely undermine the predictive accuracy of high-order PDE models.

Data-driven techniques, including machine learning, offer a promising alternative for discovering PDEs that govern such complex systems. These methods aim to automate the extraction of underlying physical laws from observational data, thus improving both the efficiency and accuracy of scientific model discovery. The goal of PDE discovery is to identify generalized, interpretable models that capture key terms and explain experimental data, while remaining robust to noise and outliers. However, several critical challenges remain: i) incomplete identification of essential variables in multiscale systems, particularly those with indistinct components, and ii) error propagation from noisy or sparse measurements, which weakens model robustness. Current methods primarily include sparse regression, symbolic regression, and hybrid approaches that incorporate deep learning architectures.

As an important tool in PDE discovery, sparse regression has been applied in many fields, especially in the automatic discovery of differential equations and system coordinates from data \cite{2,3,4}. Derived from the sparse identification of nonlinear dynamic systems (SINDy) approach, the functional identification of nonlinear dynamics (PDE-FIND) \cite{2} has brought about a significant breakthrough in the field of PDE discovery. This methodology constructs an overcomplete candidate function library and employs a sequential thresholded least squares algorithm to discover linear combinations of function terms from measurements on a regular grid. Several extensions to SINDy have improved its applicability and robustness. The SINDy with control method, incorporates control inputs to account for external drive or control signals, extending its use from uncontrolled systems to practical applications such as engineering control systems and robotics \cite{5}. Furthermore, the Lagrangian-SINDy approach characterizes system dynamics in Lagrangian form, making it particularly effective for modeling physical systems with constraints or energy conservation \cite{6}. Implicit-SINDy further extends the capability of the method by enabling the recognition of dynamical behaviors that are difficult to represent explicitly, thus improving its adaptability to complex nonlinear dynamics \cite{7}. The inclusion of weak and integral forms allows SINDy to handle noisy or incomplete data more effectively, demonstrating improved stability in the modeling of PDEs \cite{8,9}. Despite these advancements, the methods still strongly depend on a priori knowledge for constructing candidate function libraries, and the sensitivity of higher-order derivative calculations to noise limits their application in complex systems.

The concept of symbolic regression (SR) was first introduced by Koza \cite{10} as a tool for data-driven scientific discovery, using genetic programming (GP) to explore mathematical models in symbolic form. GP enables efficient search and optimization in large model spaces to identify the optimal mathematical expression by simulating the selection and optimization processes observed in biological evolution \cite{11}. 
Bongard and Lipson applied SR to discover ODE models \cite{12}, proposing a methodology that could identify dynamical equations with explanatory power, thus improving our understanding of system behavior. A further extension of SR, the symbolic genetic algorithm (SGA), combines a genetic algorithm with symbolic representations, utilizing binary trees to represent PDEs. This approach increases the flexibility and adaptability of SR for modeling complex systems \cite{13}. With advances in automation and machine learning techniques, Uderscu and Tegmark \cite{14} introduced the AI-Feynman method, which successfully identified 100 classical Feynman equations in physics by leveraging simplifying properties such as symmetry, separability, and compositeness. This demonstrated that symbolic regression can effectively uncover physically meaningful invariants and laws. Additional methods, such as AI-Aristotle and AI-Lorenz, employ symbolic regressors such as PySR to capture dynamics learned from sparse and noisy data via black-box machine learning models \cite{15,16}. Although symbolic regression has made significant strides in expressiveness and flexibility, it still faces challenges when handling highly noisy and limited data. For example, symbolic regression can generate models with redundant function terms and significant coefficient errors, while crossover and mutation operations in symbolic genetic algorithms can trigger instability in the resulting equations. Nevertheless, SR continues to show great promise in discovering invariants and simplifying expressions, particularly in physical systems with inherent symmetry and structural properties.

The application of deep learning models to dynamical systems has seen remarkable progress, particularly in dealing with complex, high-dimensional systems. Physics-Informed Neural Networks (PINNs) have demonstrated strong capabilities in simulating complex systems at low computational cost while maintaining physical constraints \cite{17,18,19}. Moreover, deep models have been employed to approximate evolutionary operators through residual networks and Eulerian approximations, with successful applications in fields such as parametric PDEs, molecular dynamics, and turbulence modeling \cite{20,21,22}. Recently, several hybrid deep learning approaches combining sparse regression, symbolic regression, and other methods have been proposed to address a range of challenges. For example, Both et al. \cite{23} introduced a hybrid method that integrates sparse regression and neural networks to recover the state of the system from noisy data. Xu et al. \cite{24} developed a deep learning genetic algorithm (DLGA-PDE) to discover parametric PDEs from sparse and noisy data, while Sun et al. \cite{25} proposed a similar deep learning genetic algorithm for constrained learning of the time evolution of the system. In molecular dynamics, Mardt et al. \cite{26} used deep learning models for scalar process modeling, and Kim et al. \cite{27} employed convolutional neural networks and finite difference methods to compute partial derivatives for time-dependent PDEs. Furthermore, Long et al. \cite{28} proposed Symbolic Regression Neural Networks (SymNet), which enhance the expressive power of models and improve the descriptive performance of complex systems by incorporating diverse activation functions.

Symbolic regression methods, when combined with deep learning models, have demonstrated excellent performance in modeling complex systems. These methods not only automate the identification of system dynamics, but also offer robust generalization capabilities. The symbolic regressor PySR \cite{29} has significantly improved the efficiency of equation discovery, drastically reducing the time required by combining genetic programming with symbolic regression techniques. To address the challenges mentioned above, this paper proposes a novel approach to PDE discovery that integrates symbolic regression and deep learning. Specifically, the method employs a multi-residual attention network to learn the response function \( u \) and its higher-order partial derivatives of the system, while utilizing the Savitzky-Golay filter \cite{30} to denoise the noisy data and enhance the robustness of the model. The proposed PySR symbolic regressor efficiently searches for dynamical equations, discovering physically meaningful explicit expressions. Compared to existing methods, our approach offers significant advantages in the following areas:

(i) \textbf{Enhanced Feature Learning Capability}. The residual attention mechanism effectively alleviates the impact of data sparsity and noise on the calculation of higher-order derivatives, improving the stability and accuracy of the model. It can effectively extract critical features and minimize the interference caused by noise, especially in complex systems.
   
(ii) \textbf{Concise and Interpretable Model Expressions}. Symbolic regression optimizes mathematical expressions by reducing redundancy and significantly minimizing coefficient errors. Compared to traditional methods, symbolic regression simplifies the expression, providing more interpretable equations that enhance model transparency.

(iii) \textbf{Efficient Computation and Robustness}. By combining the feature extraction capabilities of deep learning and the optimization advantages of symbolic regression, the proposed method demonstrates excellent performance and computational efficiency, especially when dealing with complex systems and noisy data. Specifically, compared to traditional methods (e.g., DLGA, R-DLGA \cite{31}, and NN-GA), our method outperforms in terms of both accuracy and computational efficiency. For instance, the ANN-PYSR method, based on symbolic regression, significantly reduces the relative error (0.0015) and computational time (47.86 seconds), achieving a nearly 98\% reduction compared to DLGA (1973.32 seconds).

Experimental validation on one- and two-dimensional equations shows that the proposed method achieves superior accuracy and stability even under varying noise levels. The paper is organized as follows: Section 2 introduces the ANN-PYSR framework, detailing its problem setup, data assumptions, and integrated approach combining Savitzky-Golay denoising, residual attention networks for derivatives, and PySR for equation identification. Section 3 validates the method through numerical experiments on Burgers, Chaffee-Infante, KdV, and Wave equations, evaluating relative error and computation time while comparing performance against DLGA/R-DLGA methods. Section 4 analyzes algorithm robustness under noisy/sparse sampling conditions, examining error heatmaps, sampling density (recommending ≥50\% sampling rate), and noise tolerance (up to 200\% on Burgers' equation). Finally, Section 5 summarizes ANN-PYSR's advantages in efficiency/accuracy and outlines future directions including high-dimensional extensions and adaptive sampling strategies.
\section{Frameworks}

\subsection{Problem Setting}

The primary objective of this research is to automatically discover the unknown PDEs governing the evolution of complex systems from discrete spatio-temporal observational data. This is achieved through an integrated pipeline combining advanced signal processing, deep learning-based derivative estimation, and symbolic regression. We consider a physical or engineering system whose state is described by an $m$-dimensional vector function $\bm{u}(\bm{x}, t) = [u_1(\bm{x}, t), u_2(\bm{x}, t), \dots, u_m(\bm{x}, t)]^T$. This state vector $\bm{u}$ is a function of $d$-dimensional spatial coordinates $\bm{x} = [x_1, x_2, \dots, x_d]^T \in \Omega \subset \mathbb{R}^d$ and time $t \in [0, T]$. We hypothesize that the system's dynamics are governed by a PDE system of the general form:
\begin{equation}
\frac{\partial \bm{u}(\bm{x}, t)}{\partial t} = \mathcal{N}(\bm{u}, D\bm{u}, D^2\bm{u}, \dots, D^k\bm{u}; \bm{x})
\label{eq:pde_system_general}
\end{equation}
where $\frac{\partial \bm{u}}{\partial t} = \left[ \frac{\partial u_1}{\partial t}, \dots, \frac{\partial u_m}{\partial t} \right]^T$ is the temporal derivative of the state vector. The term $D^j\bm{u}$ represents the collection of all partial derivatives of $\bm{u}$ with respect to the spatial variables $\bm{x}$ up to the order $j$. The operator $\mathcal{N} = [\mathcal{N}_1, \mathcal{N}_2, \dots, \mathcal{N}_m]^T$ is an unknown, typically nonlinear, vector-valued functional operator.

The input data available for this study consists of state observations, denoted as $\mathbf{U}_{\text{data}} = \{ \bm{u}(\bm{x}_i, t_j) \}$, at $N = N_s \times N_t$ discrete spatio-temporal sampling points $\{(\bm{x}_i, t_j) : i=1, \dots, N_s; j=1, \dots, N_t\}$. Recognizing that real-world observational data are often corrupted by noise, the initial stage of our methodology involves data preprocessing and derivative estimation. Firstly, a Savitzky-Golay (S-G) filter is applied to the raw observational data $\mathbf{U}_{\text{data}}$ to smooth and denoize, resulting in $\tilde{\mathbf{U}}_{\text{data}}$. Subsequently, to accurately estimate the temporal derivatives $\hat{\dot{\mathbf{U}}}_{ij}$ and the necessary spatial derivatives (e.g., $\hat{D\bm{u}}_{ij}, \dots, \hat{D^k\bm{u}}_{ij}$), we employ a Residual Attention Network. This network is specifically designed to learn and compute high-precision partial derivative values from (potentially sparse) spatio-temporal data. The network takes the denoised data points $\tilde{\bm{u}}(\bm{x}_i, t_j)$ and their neighborhood information as input, outputting estimates of their derivatives at each sampling point. These derivative values, obtained through S-G filtering and the Residual Attention Network, form the foundational input for the subsequent symbolic regression algorithm.

To discover the explicit mathematical expressions for each component $\mathcal{N}_l$ of the unknown operator $\mathcal{N}$, this research utilizes PySR, an efficient open-source toolkit for symbolic regression. PySR, grounded in genetic programming principles and enhanced with modern optimization techniques, automatically searches for and constructs mathematical models from a predefined set of basic building blocks (primitives), denoted as $\mathcal{P}$. This set $\mathcal{P}$ typically comprises a Variable Set $\mathcal{V}$, a Mathematical Operator Set $\mathcal{O}$ and a Constant Set $\mathcal{C}$. PySR conducts an efficient, often distributed, evolutionary search within the vast function space $\mathbb{F}$ generated by $\mathcal{P}$. It iteratively combines these building blocks to generate and evolve a multitude of candidate mathematical expressions $\mathcal{F}_{c,l}(\mathcal{V}_{ij}; \bm{c}_l)$. Here, $\mathcal{V}_{ij}$ represents the available variable values at point $(\bm{x}_i, t_j)$, and $\bm{c}_l$ is a set of constant parameters within the expression that require optimization. For each state component $u_l$, PySR aims to identify a functional form $\mathcal{F}_{c,l}$ such that $\hat{\dot{U}}_{l,ij} \approx \mathcal{F}_{c,l}(\mathcal{V}_{ij}; \bm{c}_l)$.

The core driving force of the PySR symbolic regression process is the optimization of a predefined objective function (or fitness function). This function evaluates the quality of each candidate expression $\mathcal{F}_{c,l}$, seeking a Pareto-optimal balance between the model's predictive accuracy and its simplicity (complexity). Specifically, for the $l$-th state component $u_l$, we endeavor to minimize a loss function $L(\mathcal{F}_{c,l}, \bm{c}_l)$ of the form:
\begin{equation}
L(\mathcal{F}_{c,l}, \bm{c}_l) = \underbrace{ \frac{1}{N_s N_t} \sum_{i=1}^{N_s} \sum_{j=1}^{N_t} \left( \hat{\dot{U}}_{l,ij} - \mathcal{F}_{c,l}(\mathcal{V}_{ij}; \bm{c}_l) \right)^2 }_{\text{Mean Squared Error (MSE)}} + \lambda \cdot \underbrace{ \text{Complexity}(\mathcal{F}_{c,l}) }_{\text{Model Complexity}}.
\label{eq:loss_function_pysr}
\end{equation}
In Equation~\ref{eq:loss_function_pysr}, the first term is the Mean Squared Error (MSE), which quantifies the discrepancy between the temporal derivatives predicted by the candidate expression $\mathcal{F}_{c,l}$ (given constants $\bm{c}_l$) and the temporal derivatives $\hat{\dot{U}}_{l,ij}$ estimated by the Residual Attention Network. PySR incorporates an internal mechanism for optimizing the constants $\bm{c}_l$ (i.e., finding $\bm{c}_l^* = \arg\min_{\bm{c}_l} \text{MSE}(\mathcal{F}_{c,l}, \bm{c}_l)$), typically employing numerical optimization algorithms such as L-BFGS. The second term, $\text{Complexity}(\mathcal{F}_{c,l})$, penalizes the structural complexity of the candidate expression $\mathcal{F}_{c,l}$. PySR allows for flexible definitions of this term, often based on the number of nodes in the expression tree or specific combinations of operators. The parameter $\lambda \ge 0$ is a regularization coefficient that balances the trade-off between goodness-of-fit and model parsimony. Through its evolutionary algorithm, PySR systematically searches the space of mathematical expressions to identify one or a set of Pareto-optimal expressions $\mathcal{F}_{c,l}^*$ that minimize the objective function $L(\mathcal{F}_{c,l}, \bm{c}_l^*)$.

Ultimately, for each state component $u_l$, the optimal expression $\mathcal{F}_{c,l}^*(\mathcal{V}; \bm{c}_l^*)$ identified by PySR symbolic regression constitutes an explicit mathematical approximation for the corresponding component $\mathcal{N}_l$ in the discovered PDE system. 
This specific workflow, integrating S-G filtering, deep learning-based derivative estimation, and PySR symbolic regression, enables the robust extraction of dynamical information from noisy data and facilitates the autonomous exploration and construction of mathematical structures that accurately describe and explain the observed data.

In practical scenarios, the observed data \( \tilde{u} \) is typically disturbed by noise, which is modeled as
\begin{equation}
\tilde{u} = u + \text{noise} \cdot \text{std}(u) \cdot \mathcal{N}(0, 1),
\end{equation}
where \( \text{std}(u) \) represents the standard deviation of \( u \), and \( \mathcal{N}(0, 1) \) denotes a standard normal random variable. The factor \( \text{noise} \) controls the amplitude of the noise, while \( \mathcal{N}(0, 1) \) generates a matrix of Gaussian random numbers with dimensions matching the shape of \( u \) (i.e., \( u.shape[0] \times u.shape[1] \)). This formulation allows the noise term to be scaled by both the standard deviation of \( u \) and the noise factor, thus simulating measurement errors in the data.

\subsection{Savitzky-Golay Filter Denoising Process}

In practical applications, measurement data are often subject to interference from various sources of noise, which may arise from sensor instability or external environmental factors. To accurately extract signal information, noise reduction is crucial, and the use of filters for data denoising is a common and effective technique to achieve this goal.

Consider a noisy higher-dimensional spatiotemporal dataset \( \mathbf{u}(\mathbf{x}, \mathbf{t}) \), which represents an observed signal on a spatiotemporal grid. The objective of this section is to apply the Savitzky-Golay filter for denoising the signal, while preserving as many detailed features as possible. The procedure for denoising using the Savitzky-Golay filter is outlined below \cite{32}. First, we assume that the noisy signal \( \mathbf{u}(\mathbf{x}, \mathbf{t}) \) is represented as data on a spatio-temporal grid. To facilitate the denoising process, we convert these data into a zero mean array and organize them into a flat image block matrix \( S_p(\tilde{\mathbf{u}}) \in \mathbb{R}^{p \times f} \), where \( p \) is the block size and \( f \) is the number of features. This transformation facilitates further processing using dictionary learning and sparse coding.

\textbf{Step 1: Initial signal smoothing.} Before applying the Savitzky-Golay filter, an initial smoothing of the input signal is usually required to reduce the effect of local fluctuations on the fitting results. We can use a Gaussian convolution filter to preprocess the signal, whose convolution kernel \( G(\mathbf{x}, \mathbf{t}; \sigma) \) is defined as
\begin{equation}
G(\mathbf{x}, \mathbf{t}; \sigma) = \frac{1}{(2\pi \sigma^2)^{d/2}} \exp\left( -\frac{\|\mathbf{x}\|^2 + \|\mathbf{t}\|^2}{2\sigma^2} \right),
\end{equation}
where \( \mathbf{x} \) and \( \mathbf{t} \) are vectors representing spatial and temporal coordinates, and \( d \) is the dimension of the data. The signal \( \tilde{\mathbf{U}}(\mathbf{x}, \mathbf{t}) \) is then convolved with this kernel to produce a smoothed signal \( GB(\tilde{\mathbf{U}}) \):
\begin{equation}
GB(\tilde{\mathbf{U}})(\mathbf{x}, \mathbf{t}) = \int_{-\infty}^{\infty} \int_{-\infty}^{\infty} G(\mathbf{x} - \mathbf{x}', \mathbf{t} - \mathbf{t}'; \sigma) \tilde{\mathbf{U}}(\mathbf{x}', \mathbf{t}') \, d\mathbf{x}' \, d\mathbf{t}'.
\end{equation}
The Gaussian-filtered signal serves as the basis for parameter estimation, followed by the application of the Savitzky-Golay filter for further smoothing.

\textbf{Step 2: Application of the Savitzky-Golay filter.} For additional denoising, the Savitzky-Golay filter is employed by fitting a local polynomial to the data within a sliding window. The sliding window ensures that local features are preserved while reducing noise. Assuming the window length is \( l \) and the polynomial order is \( n \), the expression of local polynomial fitting is given by:
\begin{equation}
y(\mathbf{x}, \mathbf{t}) = a_0 + a_1 \mathbf{x} + a_2 \mathbf{x}^2 + \cdots + a_n \mathbf{x}^n + b_1 \mathbf{t} + b_2 \mathbf{t}^2 + \cdots + b_n \mathbf{t}^n,
\end{equation}
where \( \{a_i\}_{i=0}^n \) and \( \{b_i\}_{i=0}^n \) are the polynomial coefficients to be solved. To minimize the fitting error, we introduce the error function \( E \), defined as the sum of the squared errors between the fitting polynomial \( y(\mathbf{x}_i, \mathbf{t}_j) \) and the observed data \( \tilde{\mathbf{U}}(\mathbf{x}_i, \mathbf{t}_j) \):
\begin{equation}
E = \sum_{i=-m}^{m} \sum_{j=-m}^{m} \left( \tilde{\mathbf{U}}(\mathbf{x}_i, \mathbf{t}_j) - y(\mathbf{x}_i, \mathbf{t}_j) \right)^2.
\end{equation}
Here, \( l = 2m + 1 \) represents the window length (odd), with \( m \) being half the width of the window. The window size ensures that there is a center for the polynomial fitting, and this is typically odd for symmetric smoothing. By performing least-squares minimization, the error function is expressed in matrix form as:
\begin{equation}
E = ||A \mathbf{a} - \tilde{\mathbf{U}}||^2,
\end{equation}
where \( A \) is the design matrix containing the polynomial basis functions for each data point, \( \mathbf{a} \) is the vector of polynomial coefficients to be determinated, and \( \tilde{\mathbf{U}} \) is the vector of observed data. Solving the normal equation:
\begin{equation}
A^T A \mathbf{a} = A^T \tilde{\mathbf{U}},
\end{equation}
yields the optimal polynomial coefficients \( \mathbf{a} \).

\textbf{Step 3: Computing Derivatives of Denoised Signals.} Once the polynomial fitting is complete, the Savitzky-Golay filter not only smooths the signal but can also be used to estimate the derivatives. The derivative of the fitted polynomial can be expressed as
\begin{equation}
\mathbf{U}'(\mathbf{x}, \mathbf{t}) = a_1 + 2a_2 \mathbf{x} + 3a_3 \mathbf{x}^2 + \cdots + n a_n \mathbf{x}^{n-1} + b_1 + 2b_2 \mathbf{t} + 3b_3 \mathbf{t}^2 + \cdots + n b_n \mathbf{t}^{n-1}.
\end{equation}
This method efficiently computes the derivatives of the signal while mitigating the effects of noise. The Savitzky-Golay filters are widely used for both smoothing and derivative estimation in noisy signals, offering excellent performance in maintaining important features while reducing noise. The ability to compute derivatives is particularly useful when analyzing the rate of change of physical systems or other dynamic signals.

\textbf{Step 4: Optimization and Selection of Denoising Parameters.} The effectiveness of the Savitzky-Golay filter depends on the appropriate selection of window length \( l \) and polynomial order \( n \). A large window length leads to significant smoothing but may obscure finer signal details, whereas a higher polynomial order can overfit the data, increasing the error. To optimize these parameters, we use the MSE to identify the optimal values of \( (\sigma^*, l^*) \):
\begin{equation}
\{\sigma^*, l^*\} = \arg \min_{\sigma, l} \frac{1}{N} \sum_{i=1}^{N} \left[ SG(\tilde{\mathbf{U}}, \sigma, l) - GB(\tilde{\mathbf{U}}) \right]^2,
\end{equation}
where \( SG(\tilde{\mathbf{U}}, \sigma, l) \) represents the signal after smoothing with the Savitzky-Golay filter, and \( GB(\tilde{\mathbf{U}}) \) is the reference signal after Gaussian filtering. By minimizing the mean square error, the optimal smoothing effect can be achieved. As shown in Figure \ref{fig1a}, the smooth solution obtained by applying the Savitzky-Golay (S-G) filter to denoise the noisy data of the KdV equation clearly demonstrates the effective restoration of the dynamic structure during the denoising process.

\subsection{Residual Attention Neural Network}

Fully connected neural networks (FCNNs) have been extensively applied in various scientific computing and machine learning tasks. However, as the depth of the network increases, the vanishing or exploding gradient problem becomes a significant bottleneck during training, severely hindering optimization. Additionally, traditional feedforward networks typically employ a fixed weight structure, which limits their capacity to capture dynamic relationships among input features. To address these challenges, we propose the residual Attention Neural Network (ANN), which incorporates two transformer networks to map input data into a high-dimensional feature space. In high-dimensional tasks such as NLP, computer vision, and scientific computing, the dual transformer networks \( U \) and \( V \) extract complementary features, enhancing the model’s adaptability to heterogeneous data and cross-modal tasks. This framework enhances the transfer of information between layers through residual connections, while simultaneously applying weighted fusion of static feature encoding for more effective integration.

The forward propagation rule introduced here combines residual paths with dynamic feature interactions, thereby significantly improving both the model's convergence speed and predictive performance. Experimental results demonstrate that, compared to conventional neural networks, the proposed architecture offers substantial improvements in fitting accuracy and convergence efficiency.

Let \( X \) represent the input data, which is mapped to two static feature spaces, \( U \) and \( V \), by two transformer networks as follows:
\begin{equation}
U = \phi(X W^1 + b^1), \quad V = \phi(X W^2 + b^2),
\end{equation}
where \( W^1, W^2 \) and \( b^1, b^2 \) are the weights and biases of the transformer networks, respectively, and \( \phi \) is a nonlinear activation function (in this case, the sine function is chosen). 

\begin{figure}[tp]
    \centering
    \begin{subfigure}[b]{0.3\textwidth}
        \centering
        \includegraphics[width=\textwidth]{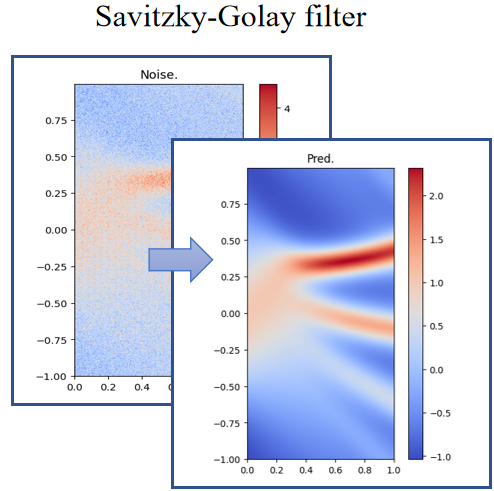} 
        \caption{} 
        \label{fig1a} 
    \end{subfigure}
    \hspace{0.035\textwidth}
    \begin{subfigure}[b]{0.65\textwidth}
        \centering
        \includegraphics[width=\textwidth]{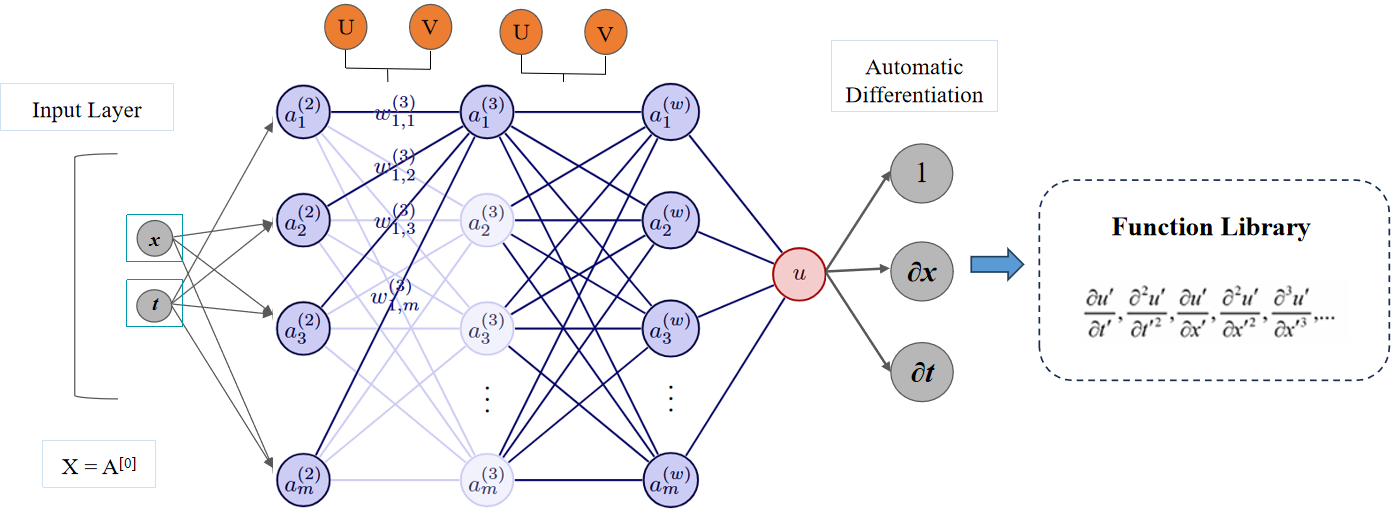} 
        \caption{} 
        \label{fig1b} 
    \end{subfigure}
    
    \vspace{0.05\textwidth} 
    
    \begin{subfigure}[b]{0.8\textwidth}
        \centering
        \includegraphics[width=\textwidth]{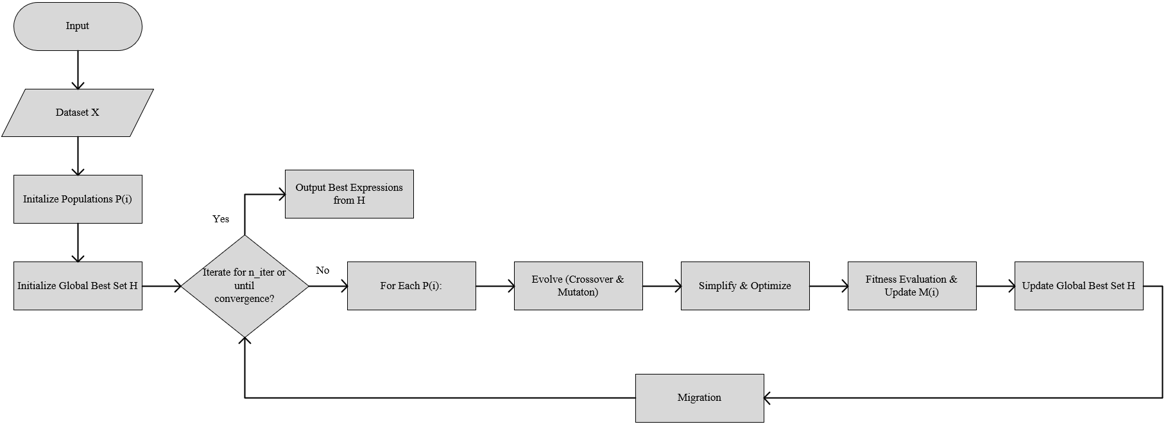} 
        \caption{} 
        \label{fig1c} 
    \end{subfigure}
    
    \caption{Overview of the algorithmic framework.  (A) Noise removal using the Savitzky-Golay filter, showing noisy data on the left and the filtered result on the right.  (B) Structure of the residual attention mechanism network, including input, hidden layers, and \( U \), \( V \) matrices in the Transformer network.  (C) Flow of the PySR symbolic regression optimization.}
    \label{fig1} 
\end{figure}

In the hidden layers, residual connections are combined with a dynamic attention mechanism. First, the previous hidden state \( H^{(k)} \) is used to calculate the dynamic attention weight \( Z^{(k)} \):
\begin{equation}
Z^{(k)} = \phi(H^{(k)} W^{z,k} + b^{z,k}), \quad k = 1, \dots, L.
\end{equation}
Here, \( W^{z,k} \) and \( b^{z,k} \) are learnable parameters. The output of the hidden layer is then updated by residual paths, dynamically fused with static features:
\begin{equation}
H^{(k+1)} = H^{(k)} + (1 - Z^{(k)}) \odot U + Z^{(k)} \odot V, \quad k = 1, \dots, L.
\end{equation}
In this equation, \( H^{(k)} \) is the hidden state of the \( k \)-th layer, \( U \) and \( V \) are the static features, \( Z^{(k)} \) is the dynamic attention weight, and \( \odot \) represents element-wise multiplication. The final output of the network is computed through a linear transformation:
\begin{equation}
u_\theta(X) = H^{(L+1)} W + b,
\end{equation}
where \( W \) and \( b \) are the weights and biases of the output layer, respectively, and \( H^{(L+1)} \) is the final hidden state. The architecture of this network is illustrated in Figure \ref{fig1c}.

\subsection{Symbolic Regression}

PySR is a multi-population evolutionary algorithm specifically designed for symbolic regression (SR), simultaneously evolving multiple independent populations to address complex regression tasks. Each population represents a set of candidate solutions, with individuals within these populations progressively optimized through evolutionary mechanisms. Extending classical evolutionary algorithms, PySR incorporates several modifications, most notably asynchronous parallelism. In this framework, each population independently undergoes evolutionary processes without requiring global synchronization, significantly enhancing the efficiency and effectiveness of search and optimization procedures. Let \(P_t^i\) denote the state of the \(i\)-th population at time step \(t\), where \(i\) indexes the populations. The evolutionary process evaluates each population’s strengths and weaknesses using a fitness function \(f(P_t^i)\). Typical fitness measures include data-fitting criteria such as MSE, with the goal of optimizing populations over successive iterations by either maximizing or minimizing this fitness function：
\begin{equation}
P_{t+1}^i = \text{evolve}(P_t^i, f),
\end{equation}
where \(\text{evolve}\) represents the evolutionary process applied to \(P_t^i\) to yield the next generation \(P_{t+1}^i\) through selection, crossover, mutation, and other operations. Tournament selection, first introduced by Goldberg et al. \cite{33} and later generalized \cite{34}, is used to select individuals based on fitness. A tournament set \(\{p_1, p_2, \dots, p_k\}\) is formed by randomly selecting \(k\) individuals from the current population, and the individual with the highest fitness is chosen as the winner:
\begin{equation}
p_{\text{winner}} = \arg \max_{p \in \{p_1, p_2, \dots, p_k\}} f(p),
\end{equation}
where \(f(p)\) represents the fitness function, which evaluates the merit of individual \(p\). The tournament size \(k\) controls the selection pressure: the larger \(k\) favors the high-fit individuals, while the smaller \(k\) introduces more randomness into the selection process.

The core operations of the evolutionary algorithm include mutation and crossover, which introduce diversity into the population to prevent premature convergence to local optima. Mutation occurs with a probability \(P_{\text{mutation}}\), randomly altering some genes in an individual to generate a new one. Mutation can be seen as applying a perturbation function to an individual \(p\):
\begin{equation}
p' = \text{mutate}(p).
\end{equation}

Crossover, in which two parents exchange genes to create offspring, ensures the inheritance of traits from both parents while promoting diversity. Let \(p_1\) and \(p_2\) be parent individuals, their crossover produces offspring \(p'_1\) and \(p'_2\):
\begin{equation}
(p'_1, p'_2) = \text{crossover}(p_1, p_2).
\end{equation}

The crossover operation facilitates the exploration of new areas in the solution space by recombining parental traits. The optimization algorithm’s search efficiency and reconciliation are enhanced through a cross-mutation algorithm that follows a cycle of ``evolution-simplification-optimization''. In the ``evolution'' phase, new generations are created via crossover and mutation, and fitness is used to evaluate individuals. The fitness function is defined as:
\begin{equation}
LE = \text{Loss}(E, X) \cdot \exp(\lambda \cdot \text{Complexity}(E)),
\end{equation}
where \(\text{Loss}(E, X)\) measures the error in individual \(E\)'s fit to the target data \(X\), \(\text{Complexity}(E)\) denotes the complexity of individual \(E\), and \(\lambda\) is a balance parameter. The optimal individual \(E^*\) is selected based on the fitness value:
\begin{equation}
E^* = \arg \min_{E \in \{E_1, E_2, \dots, E_N\}} LE(E).
\end{equation}

Equation \(P_i = \alpha_H \cdot H + \alpha_M \cdot M_j + \epsilon\) generates the next generation of populations, where \(H\) is the historically optimal individual, \(M_j\) is the current highly adapted individual and \(\epsilon\) is a random perturbation. In the ``simplification'' phase, algebraic equivalence rules simplify the equations, reducing complex expressions to simpler forms. While simplification is infrequent, it effectively reduces model complexity and limits the search space size without sacrificing solution diversity. The ``optimization'' phase refines the fitness of individuals by adjusting parameters using classical optimization algorithms, such as L-BFGS. For each individual candidate, the optimization algorithm tunes its constants to improve the accuracy of the final solution. After several rounds of ``evolution-simplification-optimization'', the output is the set of populations with optimal fitness: \(\mathcal{H} = \{E_1^*, E_2^*, \dots, E_K^*\}\), containing the best solutions that meet the optimization objective. This design achieves an effective balance between global search and local optimization, significantly improving algorithm performance and result interpretability. The algorithm flow is summarized in Algorithm \ref{tag1}.

For final model selection, PySR adopts a scoring mechanism based on both loss and complexity. Given a set of candidate models \(\{M_1, M_2, \dots, M_k\}\), an initial screening is performed based on the loss values of the models. The top models with the smallest loss values are selected, typically a small proportion \(r\) (e.g., the top 20\%) of the models, such that:
\begin{equation}
\text{Loss}_i \leq \text{Percentile}(\text{Loss}_1, \text{Loss}_2, \dots, \text{Loss}_k, r).
\end{equation}
This step filters out models with smaller loss values, ensuring better fit. For models with similar loss values or those meeting the screening criteria, a further selection is made using a score (Score) calculated as:
\begin{equation}
\text{Score} = -\frac{\log(\text{Loss}_i / \text{Loss}_{i-1})}{\text{Complexity}_i - \text{Complexity}_{i-1}},
\end{equation}
where \(\text{Loss}_i\) represents the loss of the \(i\)-th candidate model. In each round of evaluation, the algorithm computes the loss and complexity for each candidate model. \(\text{Loss}_{i-1}\) and \(\text{Complexity}_{i-1}\) represent the loss and complexity of the previous model in the sequence. The score reflects the balance between loss and complexity, with higher scores indicating models that maintain lower complexity while achieving lower loss.

\section{Results}

In this section, we apply the enhanced PySR symbolic regression algorithm, integrated with a residual attention mechanism (ANN-PYSR), to identify partial differential equations (PDEs) from key domains such as fluid dynamics and quantum mechanics. The targeted equations include the Korteweg-de Vries (KdV) equation, the Burgers equation, the Convection-Diffusion equation, the Chaffee-Infante equation, the Wave equation, and the Klein-Gordon (KG) equation. These examples demonstrate the effectiveness of the proposed framework in deriving interpretable and concise models from sparse and noisy datasets. The dataset used in this study is based on the numerical solutions to these equations.

The PySR algorithm employed in this work is enhanced with several hyperparameters to control the evolutionary process. Specifically, the number of populations \( n_k \) is set to 40, with each population containing \( L = 1000 \) expressions. The initial complexity of the random expressions is set to 3. The substitution ratio for the global optimal set, \( \alpha_H \), is 0.05, and the substitution ratio for the combinatorial set, \( \alpha_M \), is also set to 0.05. These parameters are selected based on pre-experimentation to balance the exploration and exploitation phases of the algorithm. The network structure, as proposed, is centered on the improved transformer. Two sets of initial representations, \( U \) and \( V \), are generated via linear transformations and fed into a multi-layer Transformer module to capture complex dynamic properties. The network consists of six stacked Transformer layers, each with eight attention heads to compute dynamic attention weights \( Z^{(k)} \). In particular, the \( \text{tanh} \) activation function is used to enhance nonlinearity, ensuring that the output remains within the range \([-1, 1]\). The feedforward network consists of two fully connected layers, with the hidden layer dimension expanded to 804, while a dropout rate of 0.1 is applied to prevent overfitting.

\begin{algorithm}[htbp]
\caption{Modified PySR Algorithm}
\label{tag1}
\begin{algorithmic}[1]
\Require $X$: {The dataset used to find expressions}
\Ensure The best expressions at each complexity level
\Statex
\State $n_k \gets 40$ \Comment{Number of populations}
\State $L \gets 1000$ \Comment{Number of expressions in each population}
\State $\alpha_H \gets 0.05$ \Comment{Replacement fraction using expressions from $H$}
\State $\alpha_M \gets 0.05$ \Comment{Replacement fraction using expressions from $\bigcup_i M_i$}
\Statex
\Function{pysr}{$X$}
  \For{$i = 1$ to $n_k$} \Comment{Initialize populations}
    \State Create set $P_i$ containing $L$ random expressions of complexity 3
    \State Create empty set $M_i$ to store best expressions from $P_i$
  \EndFor
  \State Create empty set $H$ to store overall best expressions
  \For{$n = 1$ to $n_{\text{iter}}$} \Comment{Optimization loop}
    \For{$i = 1$ to $n_k$} \Comment{Parallelized over workers}
      \State $P_i \gets \text{evolve}(P_i, X)$ \Comment{Evolve, simplify, optimize}
      \For{each expression $E \in P_i$}
        \State Simplify $E$
        \State Optimize constants in $E$
        \State Store updated $E$ in $P_i$
      \EndFor
      \State $M_i \gets \text{most accurate expression in } P_i \text{ at each complexity}$
    \EndFor
    \State $H \gets \text{most accurate expression in } \bigcup_i M_i \cup H \text{ at each complexity}$
    \For{$i = 1$ to $n_k$} \Comment{Migration step}
      \For{each expression $E \in P_i$}
        \If{$\text{rand()} < \alpha_H$}
          \State Replace $E$ in $P_i$ with a random expression from $H$
        \EndIf
        \If{$\text{rand()} < \alpha_M$}
          \State Replace $E$ in $P_i$ with a random expression from $\bigcup_{j \neq i} M_j$
        \EndIf
      \EndFor
    \EndFor
  \EndFor
  \State \Return $H$
\EndFunction
\end{algorithmic}
\end{algorithm}

\subsection{Evaluation Metrics}

To assess the quality of the models discovered by the algorithm, we use several evaluation metrics:

1. Relative Error. This metric quantifies the discrepancy between the discovered PDE solution and the true PDE solution. By comparing the total of squared errors between the two solutions on a discrete grid, the error can be identified as follows:  
\begin{equation}
\text{error} = \left( \frac{\sum_{i=1}^{N_x} \sum_{j=1}^{N_t} | u(x_i, t_j) - u^*(x_i, t_j) |^2}{\sum_{i=1}^{N_x} \sum_{j=1}^{N_t} | u(x_i, t_j) |^2} \right)^{1/2} {\times 100\%}.
\end{equation}
Here, \( u(x_i, t_j) \) is the true PDE solution, \( u^*(x_i, t_j) \) is the solution of the discovered PDE, and \( N_x \) and \( N_t \) are the number of spatial and temporal grid points, respectively. The accuracy and stability of numerical methods are regularly evaluated using the relative error, which quantifies the difference between the predictive solution and the true solution. Smaller error values indicate better approximations of the true equation by the discovered PDE, reflecting the effectiveness of ANN-PYSR method.

2. Relative Coefficient Error (\( e_{\text{coef}} \)).  The coefficient error measures the discrepancy between the coefficients discovered by ANN-PYSR and the true system’s coefficients. Typically, this is expressed as the \(L_{2})\) relative error between corresponding coefficients, computed as follows:
\begin{equation}
e_{\text{coef}} = \frac{1}{N_{\text{term}}} \sum_{i=1}^{N_{\text{term}}} \left| \frac{\theta_i - \theta_i^{\text{true}}}{\theta_i^{\text{true}}} \right|,
\end{equation}
where \( N_{\text{term}} \) is the number of terms in the PDEs, and \( \theta_i \) and \( \theta_i^{\text{true}} \) are the coefficients of the \(i\)-th term in the discovered and true PDEs, respectively. This coefficient error measures the disparity between the discovered and true coefficients, providing insight into the accuracy of the model's representation of the real PDEs. The coefficient error is especially important in data-driven PDEs discovery, as it reflects how well the discovered model captures the coefficient characteristics of the real model. To validate the model's robustness, we conduct multiple experiments under different noise levels to assess the impact of noise on coefficient errors. By analyzing results from several independent experiments, we obtain a more robust statistical assessment, which improves the stability and accuracy of the discovered PDEs. Thus, this error metric is crucial for evaluating the consistency of the coefficients in the PDE discovery process.

\subsection{Burgers Equation}

The Burgers equation is widely applied in various fields, including fluid dynamics, nonlinear acoustics, gas dynamics, and traffic flow. Its standard form is given by
\begin{equation}
u_t = -u u_x + \delta u_{xx}.
\end{equation}
When the diffusion coefficient \( \delta \) is small, the Burgers equation exhibits pronounced nonlinear effects and slower diffusion dynamics. Physically, this means that the propagation and diffusion of physical quantities become slower and more concentrated, and instabilities and vortices in the system become more significant, leading to a more pronounced formation of excitations. To explore these phenomena, different values of the diffusion coefficient \( \delta \) were selected: 0.02, 0.04, 0.08 and 0.1. The Burgers equation was numerically solved and modeled for these parameter values. The dataset consists of 256 spatial and 201 temporal observations, covering a spatial range of \( x \in (-8, 8) \) and a temporal range of \( t \in (0, 10) \), resulting in a total dataset size of 51,456. For model training, the data was split into training and validation sets in an 8:2 ratio. During pre-training, noisy data were fitted using the Adam and L-BFGS optimizers, and training was terminated based on the performance on the validation set.

\begin{figure}[ht]
    \centering
    \includegraphics[width=1\linewidth]{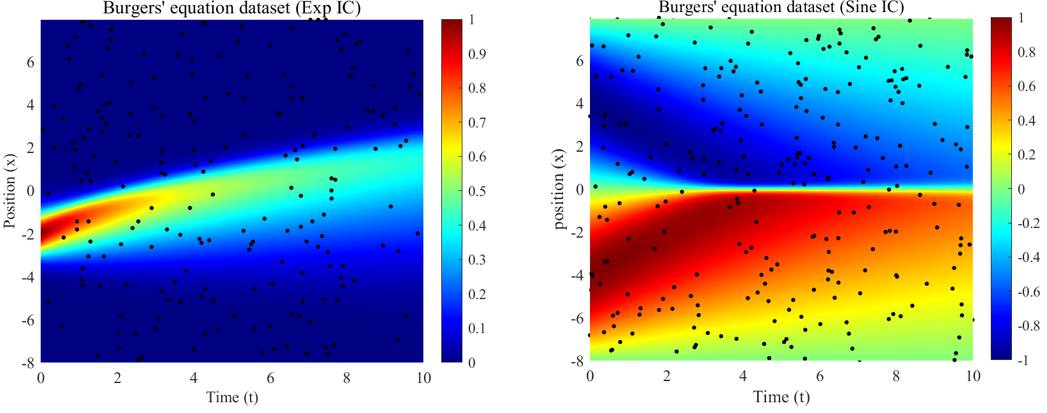}
    \caption{Burgers equation dataset.}
    \label{fig2}
\end{figure}

Figure \ref{fig2} shows the solutions for the Burgers equation dataset with \( \delta = 0.1 \), using an exponential (Exp IC) and sine (Sine IC) initial condition. The heatmap in the figure illustrates how the solution evolves over time (from 0 to 10) and space (from -8 to 8). A number of black dots, representing selected sampling points, are randomly distributed within the figure to facilitate further analysis of the solution's behavior at different time and spatial locations. This figure clearly illustrates the dynamic behavior of the Burgers equation solution. To better illustrate the stability and accuracy of the algorithm under noisy conditions, we introduced Gaussian noise ranging from 0\% to 250\% into the dataset.

\begin{figure}
    \centering
    \includegraphics[width=0.8\linewidth]{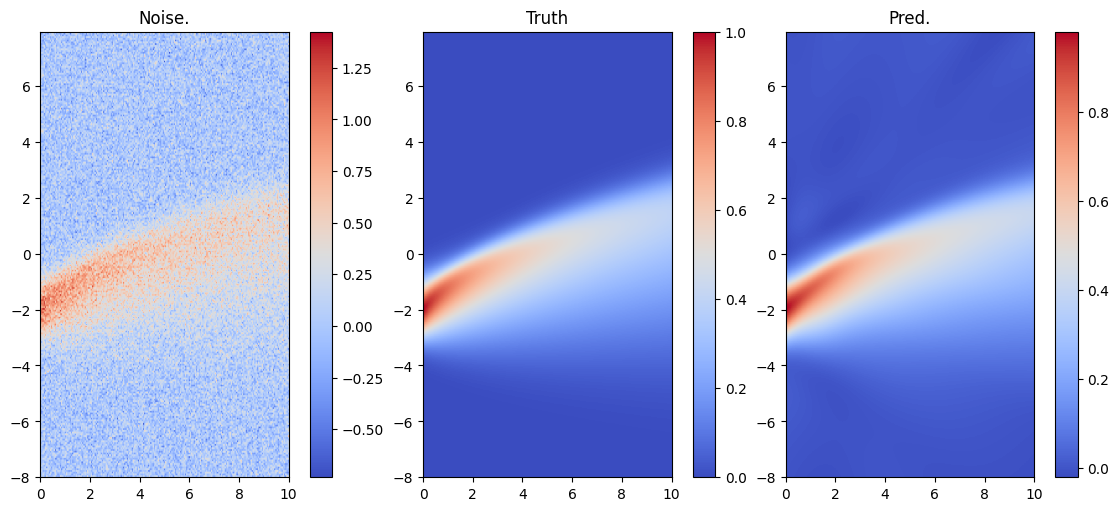}
    \caption{Comparison of solutions to the Burgers equation at a noise level of 80\%. The left side shows the input data containing noise (Noise), the center shows the true solution (Truth), and the right side shows the solution obtained by prediction through symbolic regression methods (Pred).}
    \label{fig3}
\end{figure}

We conducted systematic identification experiments on the parameterized Burgers equation using the ANN-PYSR algorithm, focusing on evaluating the algorithm’s performance across different diffusion coefficients \( \delta \) and noise levels. Figure \ref{fig3} presents the identification results for the Burgers equation (Exp IC) at a noise level of 80\%. The columns display the observed values under noisy conditions, the true solution, and the predicted solution generated by the symbolic regression method. The results demonstrate that the ANN-PYSR method successfully captures the true equation's structure, even under high noise conditions.

In Table \ref{tagg1}, we compare the performance of different candidate equations under 200\% noise for the exponential initial conditions (top) and 50\% noise for the sine initial conditions (bottom). This enables us to assess the robustness of the symbolic regression process across varying noise levels. The table presents key metrics—Complexity, Loss, and Score—for each candidate equation. The equation highlighted in bold was selected as the final model because of its optimal balance between loss and complexity, providing an accurate fit to the Burgers equation. Notably, the identification process successfully recognized the Burgers equation under both exponential and sine initial conditions, even in the presence of significant noise. The consistent identification of the equation across both noise scenarios highlights the robustness and adaptability of the PySR framework in handling noisy data, effectively capturing the Burgers equation from both initial conditions.

Table \ref{tag2} further illustrates the recognition equations at various noise levels (Exp IC), corresponding to different diffusion coefficients: \( \delta = 0.1 \), \( \delta = 0.08 \), \( \delta = 0.04 \), and \( \delta = 0.02 \). As \( \delta \) decreases, the system’s diffusion effect increases, but ANN-PYSR maintains a low parameter estimation error, even as noise influences the recognition process. For example, when \( \delta = 0.08 \), even if the noise level is 1, the parameter error of the algorithm remains below 28\%, showing excellent adaptability and stability. Moreover, at \( \delta = 0.1 \) and a noise level of 280\%, the equation structure remains recognizable within the error range. In comparison with other algorithms (such as DLGA and R-DLGA), ANN-PYSR not only performs well in recognition accuracy, but also has significantly shorter computation time, which improves the overall efficiency (see Table \ref{tag3}). Both DLGA and R-DLGA, traditional genetic programming-based PDE discovery methods, rely on genetic algorithms for equation optimization and evolution. However, these methods suffer from long computation times due to their dependence on multi-generation computation. In contrast, ANN-PYSR integrates the benefits of deep learning and symbolic regression by combining multilayer perceptron (MLP) architecture with the PySR algorithm. However, ANN-PYSR performs better than the NN-GA network architecture, which generates coarser derivative data, leading to inferior results in PDE recognition.

\begin{table}[h!]
\centering
\caption{Complexity, Loss and Score of Burgers Equations under Exponential Initial Conditions with 200\% Noise (Top) and Sine Initial Conditions with 50\% Noise (Bottom)}
\label{tagg1}
\begin{adjustbox}{max width=\textwidth}
\begin{tabular}{ccccc}
\hline
\toprule
Complexity & Loss & Score & Identified PDE \\
\midrule
1 & 0.004547 & 0.000000 & $u_t = -0.0003$ \\
3 & 0.000906 & 0.205982 & $u_t = -0.3885 \cdot u_x$ \\
5 & 0.000906 & 0.000004 & $u_t = (-0.3885 \cdot u_x) + 0.0001$ \\
7 & 0.000717 & 0.116788 & $u_t = -0.4298 \cdot u - 0.2247 \cdot u_x$ \\
\textbf{9} & \textbf{0.000331} & \textbf{0.385830} & $\mathbf{u_t = -0.9898 \cdot u_x \cdot u + 0.0981 \cdot u_{xx}}$ \\
11 & 0.000287 & 0.071967 & $u_t = 0.0863 \cdot u_{xx} - 0.6823 \cdot u \cdot u_x - 0.1140$ \\
13 & 0.000277 & 0.016679 & $u_t = 0.0863 \cdot u_{xx} - 0.1062 \cdot u_x - 0.3749 \cdot u \cdot u_x$ \\
15 & 0.000229 & 0.096307 & $u_t = -0.1528 \cdot u_{xx} \cdot u_x - 0.5477 \cdot u_x - 1.0085 \cdot u \cdot u_x$ \\
\bottomrule
\hline
\toprule
Complexity & Loss & Score & Identified PDE \\
\midrule
1 & 0.048182 & 0.000000 & $u_t = u_x$ \\
3 & 0.001051 & 0.012501 & $u_t = u_x \cdot -0.3899$ \\
5 & 0.001051 & 0.000016 & $u_t = -0.3900 \cdot u_x + 0.0002$ \\
7 & 0.000870 & 0.094488 & $u_t = -0.4138 \cdot u_x \cdot u + 0.2329 \cdot u_x$ \\
\textbf{9} & \textbf{0.000521} & \textbf{0.126116} & $\mathbf{u_t = -0.9615 \cdot u_x \cdot u + 0.0960 \cdot u_{xx}}$ \\
11 & 0.000321 & 0.000347 & $u_t = -0.9314 \cdot u \cdot u_x + 0.0548 \cdot u_{xx} \cdot u$ \\
13 & 0.000407 & 0.103315 & $u_t = -1.0207 \cdot u \cdot u_x - 0.3790 \cdot u_{xx} - 0.1540 \cdot u_x$ \\
15 & 0.000402 & 0.006198 & $u_t = -1.0207 \cdot u \cdot u_x + 1.4541 \cdot u_x \cdot u_{xx} + 0.0681 \cdot u_{xx}$ \\
\bottomrule
\hline
\end{tabular}
\end{adjustbox}
\end{table}

\begin{table}[h!]
\centering
\caption{Noise Level, Identified PDE, Time of SR and Average Coefficient Error of Burgers equation}
\label{tag2}
\begin{adjustbox}{max width=\textwidth}
\begin{tabular}{|c|c|c|c|c|}
\hline
\textbf{Equation form} & \textbf{Noise Level} & \textbf{Identified PDE} & \textbf{Time of SR (s)} & \textbf{Average Coefficient Error} \\
\hline
\multirow{10}{*}{$u_t = - u u_x + 0.1 u_{xx}$} & 0.2 & $u_t = -0.9924 u u_x + 0.1027 u_{xx}$ & 65.33 & 0.0173 \\
 & 0.5 & $u_t = -0.9796 u u_x + 0.0970 u_{xx}$ & 64.32 & 0.0252 \\
 & 0.9 & $u_t = -0.9373 u u_x + 0.0982 u_{xx}$ & 61.61 & 0.0404 \\
 & 1.2 & $u_t = -0.9769 u u_x + 0.1037 u_{xx}$ & 64.52 & 0.033 \\
 & 1.5 & $u_t = -0.9834 u u_x + 0.1030 u_{xx}$ & 63.84 & 0.0301 \\
 & 1.7 & $u_t = -0.9898 u u_x + 0.0981 u_{xx}$ & 65.84 & 0.0224 \\
 & 2.0 & $u_t = -0.9690 u u_x + 0.1013 u_{xx}$ & 62.39 & 0.0157 \\
 & 2.2 & $u_t = -0.9615 u u_x + 0.0955 u_{xx}$ & 64.07 & 0.0211 \\
 & 2.5 & $u_t = -0.9367 u u_x + 0.1018 u_{xx}$ & 63.65 & 0.0353 \\
 & 2.8 & $u_t = -0.9367 u u_x + 0.1018 u_{xx}$ & 63.65 & 0.0427 \\
\hline
\multirow{4}{*}{$u_t = - u u_x + 0.08 u_{xx}$} & 0.2 & $u_t = -0.9432 u u_x + 0.0745 u_{xx}$ & 59.85 & 0.0660 \\
 & 0.5 & $u_t = -0.9301 u u_x + 0.0732 u_{xx}$ & 62.55 & 0.0788 \\
 & 0.7 & $u_t = -0.8942 u u_x + 0.0703 u_{xx}$ & 59.33 & 0.1155 \\
 & 1.0 & $u_t = -0.7484 u u_x + 0.0552 u_{xx}$ & 62.12 & 0.2823 \\
\hline
\multirow{6}{*}{$u_t = - u u_x + 0.04 u_{xx}$} & 0.1 & $u_t = -0.9176 u u_x + 0.0356 u_{xx}$ & 61.30 & 0.0917 \\
 & 0.2 & $u_t = -0.8907 u u_x + 0.0344 u_{xx}$ & 59.86 & 0.1247 \\
 & 0.3 & $u_t = -0.8279 u u_x + 0.0381 u_{xx}$ & 64.21 & 0.1098 \\
 & 0.4 & $u_t = -0.8546 u u_x + 0.0275 u_{xx}$ & 59.30 & 0.2665 \\
 & 0.5 & $u_t = -0.7235 u u_x + 0.0266 u_{xx}$ & 58.47 & 0.3058 \\
 & 0.6 & $u_t = -0.6287 u u_x + 0.0223 u_{xx}$ & 60.25 & 0.4069 \\
\hline
\multirow{4}{*}{$u_t = - u u_x + 0.02 u_{xx}$} & 0.1 & $u_t = -0.9099 u u_x + 0.0172 u_{xx}$ & 66.68 & 0.1196 \\
 & 0.2 & $u_t = -0.8077 u u_x + 0.0149 u_{xx}$ & 64.17 & 0.2237 \\
 & 0.15 & $u_t = -0.8964 u u_x + 0.0154 u_{xx}$ & 63.99 & 0.1633 \\
 & 0.25 & $u_t = -0.7487 u u_x + 0.0161 u_{xx}$ & 61.43 & 0.2231 \\
\hline
\end{tabular}
\end{adjustbox}
\end{table}

\begin{table}
\centering
\caption{Comparison of Identified Burgers Equation (30\%noise)}
\label{tag3}
\begin{tabular}{cccc}
\toprule
Methods   & Average Coefficient Error   & Identified PDE & Time (s) \\
\midrule
DLGA      & 0.4000   & $u_t = -0.7012u u_x + 0.0509u_{xx}$ & 1973.32 \\
R-DLGA    & 0.0074   & $u_t = -0.9920u u_x + 0.0995u_{xx}$ & 2162.71 \\
NN-GA     & 0.1611   & $u_t = -0.8751u u_x + 0.0724u_{xx}$ & 40.58 \\
ANN-PYSR    & 0.0035   & $u_t = -0.9971u u_x + 0.1008u_{xx}$ & 47.86 \\
\bottomrule
\end{tabular}
\end{table}

\subsection{Korteweg-de Vries Equation}

The Korteweg-de Vries (KDV) equation is a classical nonlinear PDE that was originally introduced by Korteweg and de Vries to model the propagation of shallow water waves under weakly nonlinear and weakly dispersive conditions. The KdV equation examined in this study is
\begin{equation}
u_t = -u u_x - 0.0025 u_{xxx},
\end{equation}
where \( u(t, x) \) denotes the wave height at time \( t \) and spatial position \( x \). The nonlinear term \( -u u_x \) and the higher-order dispersion term \( -0.0025 u_{xxx} \) together govern the speed, amplitude, and stability of the wave as it propagates. Accurate identification of the parameters and structure of this classical equation is crucial for understanding the dynamical behavior of nonlinear waves, predicting the evolution of the wave field, and devising control strategies. 

\begin{figure}
    \centering
    \includegraphics[width=0.8\linewidth]{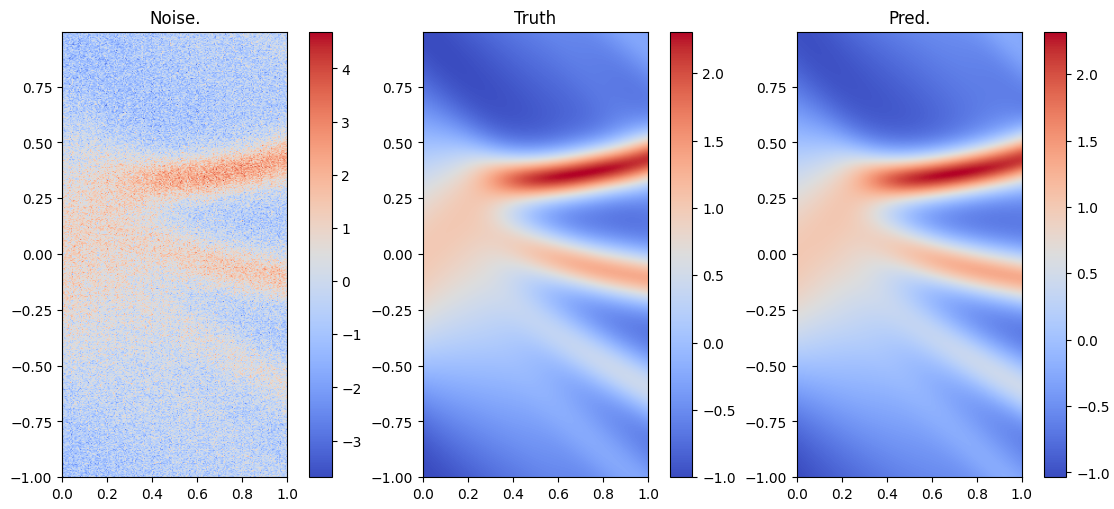}
    \caption{Comparison of the solutions of the KDV equation at a noise level of 20\%. The left side shows the input data containing noise (Noise), the center shows the true solution (Truth), and the right side shows the solution obtained by prediction through symbolic regression method (Pred).}
    \label{fig4}
\end{figure}

For this study, 512 discrete points were selected in the spatial direction, and 201 discrete points were selected in the temporal direction, resulting in a dataset of 102,912 data points as the benchmark. Subsequently, we added 0\% to 35\% Gaussian noise to assess whether the core structure and parameters of the original KDV equation could still be accurately identified under perturbed conditions. In Figure \ref{fig4}, a set of plots is presented showing the spatio-temporal dynamics of KDV equation under 20\% noise.  The columns display, in order: the noisy observations, the true solution, and the predicted solution obtained through symbolic regression. Despite the presence of noise, the method accurately recovers the underlying structure of the true equation.

Table \ref{tag4} demonstrates that the ANN-PYSR method can reliably reconstruct the core structure of the KDV equation even with noise levels as high as 20\% and 35\%. The average error of the identified parameters remains within a reasonable range of 5\% to 25\%, highlighting the algorithm's robustness under significant noise interference. Furthermore, the computational time required for ANN-PYSR shows minimal variation across different noise levels, consistently staying within an acceptable range. Specifically, symbolic regression (SR) time fluctuates slightly between 42.95s and 111.92s, indicating the method's efficiency remains stable despite increasing noise levels, further confirming its computational reliability.

Table \ref{tag5} provides a more thorough examination of the metrics for complexity, loss, and composite score under 10\% noise. These results reveal that ANN-PYSR efficiently finds candidate equations that balance accuracy and simplicity within a limited number of iterations. The bold-labeled selected results exhibit low loss errors while maintaining moderate complexity, thus avoiding overfitting and ensuring both physical interpretability and generalizability of the identified equation form. 

Comparisons with alternative methods (Table \ref{tag6}) reveal that ANN-PYSR consistently outperforms its counterparts in terms of stability. In the presence of 10\% levels, competing algorithms show noticeable degradation in recognition accuracy, while ANN-PYSR maintains its robustness. Furthermore, ANN-PYSR substantially reduces computational time, exhibiting markedly higher efficiency compared to traditional methods such as DLGA and R-DLGA, and it performs on par with or surpasses NN-GA. This superior efficiency is attributable to the optimization strategy employed by ANN-PYSR in both search space exploration and parameter estimation, facilitating rapid convergence and the identification of the joint solution within a complex parameter space. Overall, the experimental results presented in this study demonstrate that ANN-PYSR strikes an exceptional balance between accuracy, model simplicity, and computational efficiency in identifying the KDV equation, even under challenging high-noise conditions.

\begin{table}[h!] 
\centering
\caption{Noise Level, Identified PDE, Time of SR, and Average Coefficient Error of KDV equation} 
\label{tag4}
\begin{adjustbox}{max width=\textwidth} 
\begin{tabular}{cccc} 
\toprule 
Noise Level & Identified PDE & SR (time) & Average Coefficient Error \\ 
\midrule 
0.01 & $u_t = -0.9442u{u_x} - 0.002353u_{xxx}$ & 42.95s & 0.0571 \\ 
0.05 & $u_t = -0.9245u{u_x} - 0.002158u_{xxx}$ & 107.99s & 0.1060 \\ 
0.10 & $u_t = -1.0065u{u_x} - 0.002489u_{xxx}$ & 111.92s & 0.0654 \\ 
0.20 & $u_t = -0.8516u{u_x} - 0.002091u_{xxx}$ & 104.91s & 0.1559 \\ 
0.30 & $u_t = -0.7934u{u_x} - 0.002033u_{xxx}$ & 106.14s & 0.1965 \\ 
0.35 & $u_t = -0.7203u{u_x} - 0.001931u_{xxx}$ & 104.66s & 0.2535 \\ 
\bottomrule 
\end{tabular} 
\end{adjustbox} 
\end{table} 

\begin{table}[h!]
\centering
\caption{The complexity, loss, and score of the KDV equations identified under 10\% noise}
\label{tag5}
\begin{adjustbox}{max width=\textwidth}
\begin{tabular}{ccccc}
\toprule
Complexity & Loss & Score & Identified PDE \\
\midrule
1 & 3.007269 & 0.000000 & $u_t = 0.0012$ \\
3 & 2.921148 & 0.014527 & $u_t = -0.0026 \cdot u_{xx}$ \\
5 & 2.633561 & 0.051819 & $u_t = -0.0903 \cdot u_x \cdot u$ \\
7 & 2.353454 & 0.056226 & $u_t = 0.0006 \cdot -1.2681 \cdot u + 0.0006 \cdot u_{xxx}$ \\
\textbf{9} & \textbf{0.127244} & \textbf{1.458765} & $\mathbf{u_t = -1.0865 \cdot u \cdot u_x - 0.0024 \cdot u_{xxx}}$ \\
11 & 0.127145 & 0.000390 & $u_t = -0.9581 \cdot u \cdot u_x - 0.0024 \cdot u_{xxx} + 0.0100$ \\
13 & 0.126491 & 0.002576 & $u_t = 400.7806 \cdot u_x \cdot u - 1.1664 \cdot u \cdot u_x - 0.0024 \cdot u_{xxx}$ \\
15 & 0.116296 & 0.042019 & $u_t = u_{xxx} + -0.0606 \cdot u_{xx} \cdot -0.3418 \cdot u_x \cdot u \cdot -0.0024$ \\
17 & 0.112698 & 0.015714 & $u_t = u_{xxx} + 0.5001 \cdot u + -0.0606 \cdot u_{xx} \cdot -0.2384 \cdot u \cdot u_x \cdot -0.0024$ \\
\bottomrule
\end{tabular}
\end{adjustbox}
\end{table}

\begin{table}
\centering
\caption{Comparison of Identified KDV Equation (10\%noise)}
\label{tag6}
\begin{tabular}{cccc}
\toprule
Methods   & Average Coefficient Error   & Identified PDE & Time (s) \\
\midrule
DLGA      & 0.1460   & $u_t = -0.8701u{u_x} - 0.001974u_{xxx}$ & 2812.43 \\
R-DLGA    & 0.0064   & $u_t = − 0.993u{u_x} - 0.00249u_{xxx}$ & 2745.72 \\
NN-GA     & 0.2480   & $u_t = -0.7504u{u_x} - 0.001896u_{xxx}$ & 153.16 \\
ANN-PYSR    & 0.0061   & $u_t = -1.0065u{u_x} - 0.002489u_{xxx}$ & 150.29 \\
\bottomrule
\end{tabular}
\end{table}

\subsection{Klein-Gordon Equation}

The Klein-Gordon (KG) equation, which was first proposed by Oskar Klein and Walter Gordon in 1926 to describe the behavior of electrons in relativistic settings, is written as
\begin{equation}
u_{tt} = 0.5 u_{xx} - 5u.
\end{equation}
The dataset is grid data of 201 spatial observation points in the domain \( x \in [-1,1] \) and 201 temporal observation points in the domain \( t \in [0,3] \), and thus the data size is 40,401.
In the symbolic regression modeling process, for the physical properties of the Klein-Gordon equation, its proxy model is chosen to have the second-order derivative of time \( u_{tt} \) as the left end term, formalized as
\begin{equation}
u_{tt} = \Omega(u, u_x, u_{xx}, \dots).
\end{equation}
The function \( \Omega \) depends on the displacement variable \( u \), its spatial gradient \( u_x \), and the second spatial derivative \( u_{xx} \). The inclusion of \( u_{tt} \) as the primary term in the equation emphasizes the role of acceleration (i.e., the second-order time derivative) in the behavior of the scalar field, which is influenced by both its spatial distribution and internal dynamics. This formulation effectively captures the behavior of the Klein-Gordon equation as a hyperbolic PDE, describing not only the evolution of the field but also the propagation of oscillatory phenomena. The precise identification of \( \Omega(u, u_x, u_{xx}, \dots) \) through symbolic regression can offer valuable insights into the dynamic properties and underlying mechanisms of the equation, ultimately revealing the governing physical laws within a mathematical framework.

\begin{figure}
    \centering
    \includegraphics[width=0.8\linewidth]{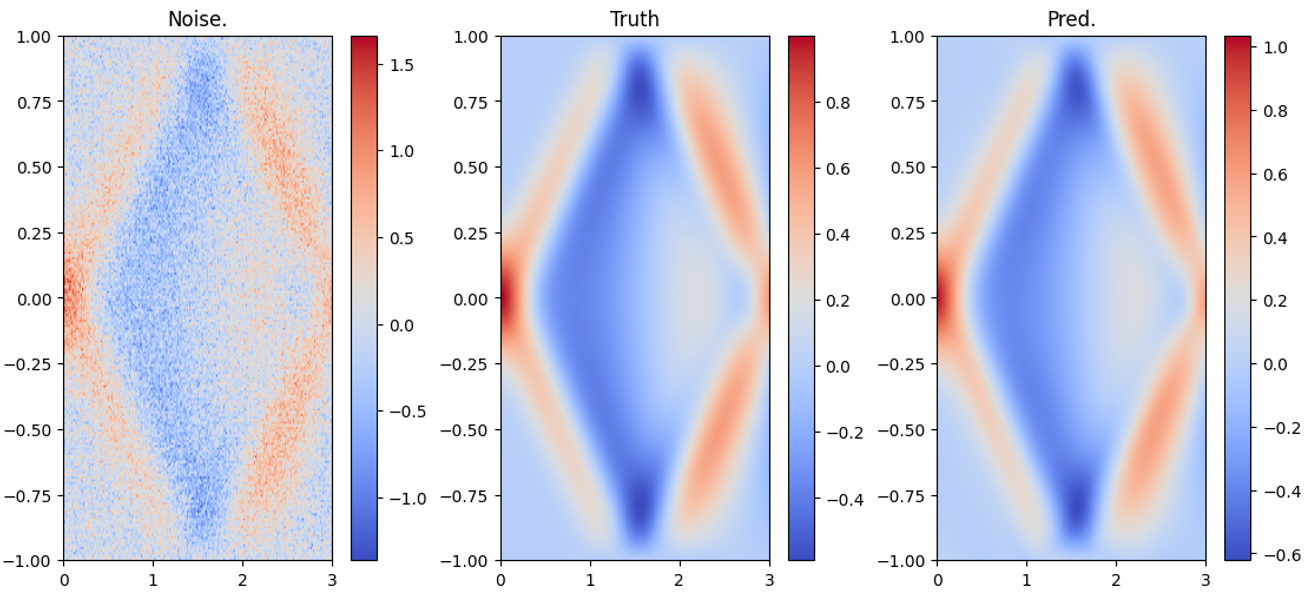}
    \caption{Comparison of the solutions of the KG equation at a noise level of 50\%. The left side shows the input data containing noise (Noise), the center shows the true solution (Truth), and the right side shows the solution obtained by prediction through symbolic regression method (Pred).}
    \label{fig5}
\end{figure}

\begin{figure}
    \centering
    \includegraphics[width=0.7\linewidth]{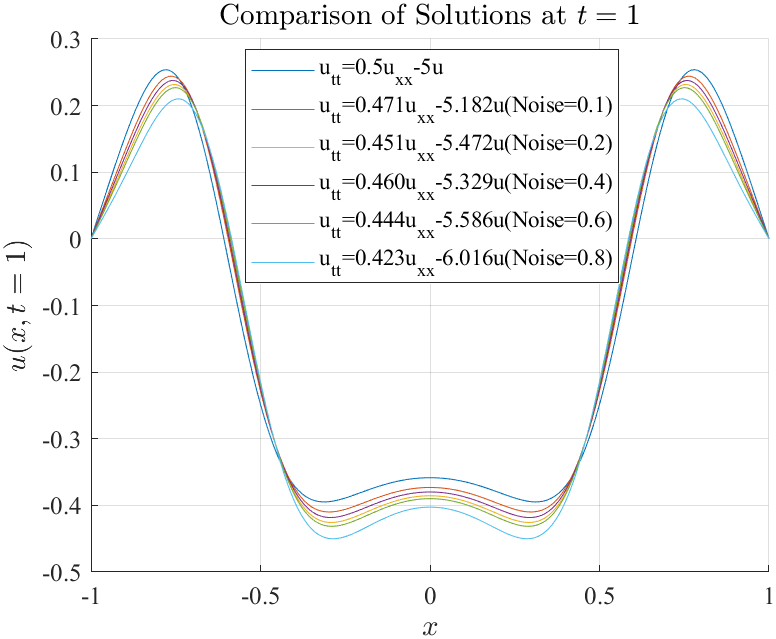}
    \caption{The plot provides a comparative analysis of the solution distributions for the KG equation recognized using the ANN-PYSR method under different noise conditions.}
    \label{fig6}
\end{figure}

In this subsection, we leverage the ANN‑PYSR algorithm to infer the symbolic form of the noisy Klein--Gordon equation, robustly recovering its governing structure under varying noise intensities. Figure~\ref{fig5} shows (from left to right) the 50\%‑noisy data, the true solution, and the symbolic regression result. Remarkably, even at elevated noise levels, our method accurately captures the system’s dynamical behavior. Figure~\ref{fig6} illustrates the spatial distribution of the identified solution compared to the true solution at time \( t = 1 \) across different noise levels. It is evident that the predicted distribution closely approximates the true one, even when \( \text{Noise} = 0.8 \). Notably, for \( \text{Noise} \leq 0.4 \), the recovered equations match the ground truth with an average coefficient error of approximately 1\%, underscoring the robustness of the ANN‑PYSR algorithm in symbolic regression for noisy data.

Furthermore, in the case of a noise level \(\ noise = 0.3\), Table~\ref{tag7} details the symbolic regression search trajectory, listing for each candidate model its complexity, training loss, and regression score. As the complexity increases, the algorithm produces more accurate models. However, when complexity becomes excessive (e.g., \( \text{Complexity} > 9 \)), the physical interpretability of the model is significantly compromised. In contrast, models with lower complexity (e.g., \( \text{Complexity} = 7 \)) not only better capture the underlying physical laws but also achieve the highest scores (\( \text{Score} = 0.4827 \)), making them the preferred choice for the final recognized equation form.

Additionally, the average runtime across the five noise levels was approximately 87.25 seconds, further confirming the efficiency of the algorithm. These results demonstrate that the ANN-PYSR algorithm strikes an excellent balance between model complexity and accuracy, efficiently generating physically meaningful and predictive models.

\begin{table}[h!]
\centering
\caption{The complexity, loss, and score of the KG equations identified under 30\% noise}
\label{tag7}
\begin{adjustbox}{max width=\textwidth}
\begin{tabular}{ccccc}
\toprule
Complexity & Loss & Score & Identified PDE \\
\midrule
1  & 17.933779 & 0.000000 & $u_{tt} = u_{xx}$ \\
3  & 2.292445 & 0.128534 & $u_{tt} = u_{xx} - 0.5692$ \\
5  & 2.267100 & 0.005559 & $u_{tt} = 0.5689 \cdot u_{xx} - 0.1592$ \\
\textbf{7}  & \textbf{0.863306} & \textbf{0.482744} & $\mathbf{u_{tt} = 0.4816 \cdot u_{xx} - 5.4343 \cdot u}$ \\
9  & 0.863083 & 0.000013 & $u_{tt} = -5.4371 \cdot u + 0.4816 \cdot u_{xx} + 0.0048$ \\
11 & 0.863893 & 0.010135 & $u_{tt} = -5.5302 + u^3 + 0.4832 \cdot u_{xx}$ \\
13 & 0.832423 & 0.008067 & $u_{tt} = -2.8771 + 2u^2 + 0.4825 \cdot u_{xx}$ \\
15 & 0.777494 & 0.034132 & $u_{tt} = 0.4721 \cdot u_{xx} - 5.8067 \cdot u + 0.1386 \cdot u_{x}^2 - 0.1676$ \\
19 & 0.764458 & 0.000120 & $u_{tt} = 0.4742 \cdot u_{xx} - 5.7646 \cdot u + 0.1413 \cdot u_{x}^2 + 0.2385 \cdot u_{x} + 0.2435$ \\
\bottomrule
\end{tabular}
\end{adjustbox}
\end{table}

\subsection{Convection-Diffusion Equation}

The convection-diffusion equations are widely utilized to model the transport of substances (e.g., pollutants, heat, or solutes) in fluid systems. These equations hold significant theoretical importance and practical applications in fields such as fluid mechanics, environmental science, and heat transfer. In this study, we consider the following form of the convection-diffusion equation:
\begin{equation}
    u_t = -u_x + 0.25u_{xx},
\end{equation}
where \(0.25\) represents the diffusion coefficient, quantifying the intensity of the diffusion effect in space. To thoroughly investigate the nature of this equation, we constructed a high-resolution gridded dataset based on a numerical solution method. The spatial domain was discretized into 256 observation points (\(x \in [0, 2]\)), and the temporal domain was divided into 100 time points (\(t \in [0, 1]\)), generating a total of 25,600 data samples. Gaussian noise with amplitudes ranging from 0 to 160\% was introduced into the original data to assess the model's ability to identify the equation under different noise levels.

\begin{figure}
    \centering
    \includegraphics[width=0.8\linewidth]{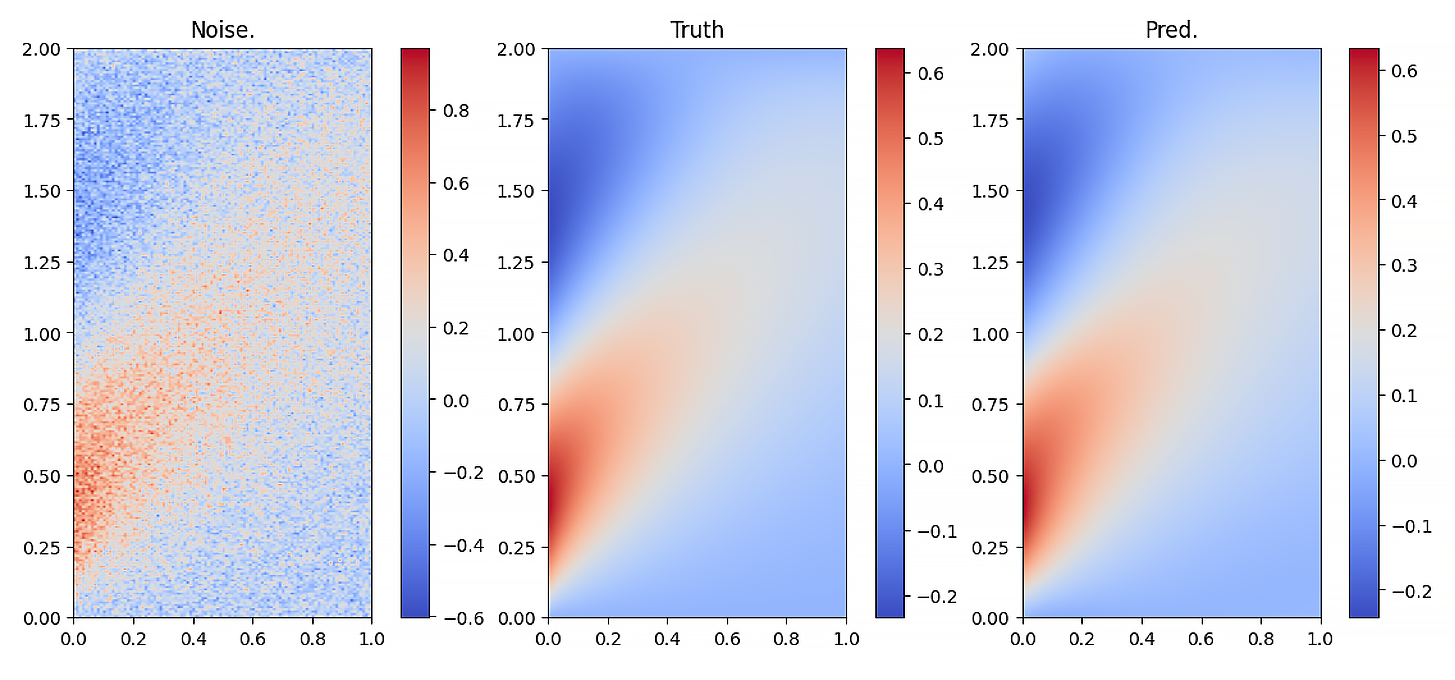}
    \caption{Comparison of the convection-diffusion equation solutions under 50\% noise. Left: Noisy input data (Noise), Center: True solution (Truth), Right: Predicted solution (Pred) obtained via symbolic regression.}
    \label{fig7}
\end{figure}

By comparing the noisy observations, true solutions, and predicted solutions (Figure \ref{fig7}), we observe that the ANN-PYSR method effectively captures the fundamental structure of the convection-diffusion equation. The predicted solutions exhibit strong agreement with the true solutions, with the \(L_2\) relative error between the predicted and true solutions reaching \(2.464 \times 10^{-3}\).

\begin{figure}
    \centering
    \includegraphics[width=0.8\linewidth]{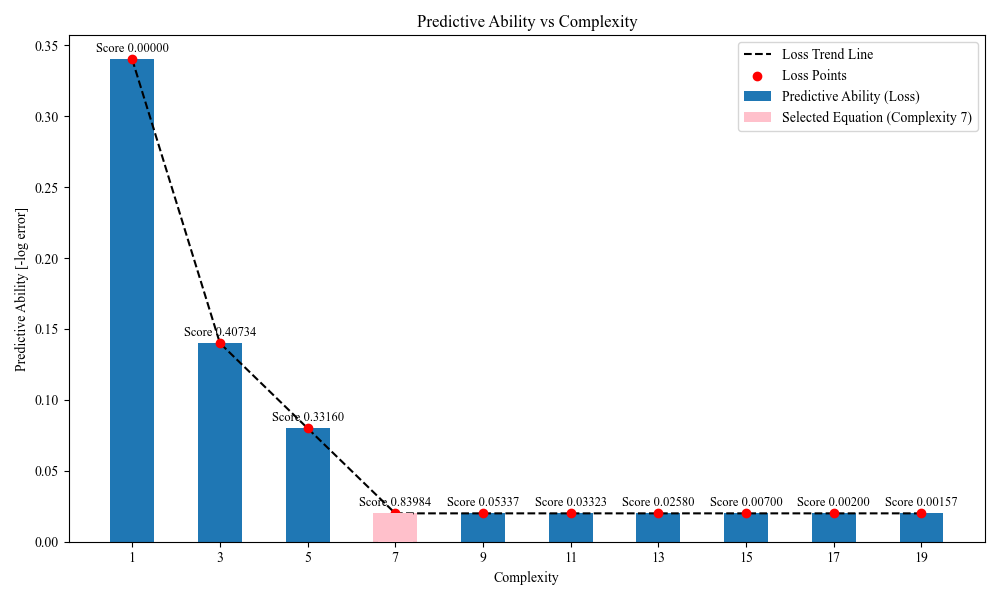}
    \caption{Identification results for the convection-diffusion equation using the ANN-PYSR algorithm, illustrating the relationship between prediction error and model complexity.}
    \label{fig8}
\end{figure}

We applied the ANN‑PYSR algorithm to automatically identify the one‑dimensional convection–diffusion equation, with results summarized in Figure~\ref{fig8} and Table~\ref{tag8}. To quantify model performance, we examined the relationship between prediction error (y axis) and model complexity (x axis). The black curve shows that, although increasing complexity yields a substantial reduction in prediction error, the magnitude of improvement progressively diminishes, indicating diminishing marginal returns. This finding suggests that model selection must balance simplicity and accuracy to avoid overfitting. Table~\ref{tag8} lists the equation form at each complexity level. Ultimately, the complexity‑7 model (highlighted in bold) achieved the highest score of 0.83984, outperforming all other candidates and striking an optimal compromise between error minimization and structural parsimony.

To further assess the robustness of the algorithm under noisy conditions, we introduced various noise levels (ranging from 0.2 to 1.6) into the dataset. Table \ref{tag9} summarizes the identified forms of the PDEs and the corresponding parameter errors at different noise levels. The results show that the core structure of the identified equations remains stable as the noise level increases, consistently following the form \( u_t = -\alpha u_x + \beta u_{xx} \). While small fluctuations in the parameters occur due to noise, the errors remain minimal. For instance, under low noise conditions (noise level 0.2), the identified equation is \( -0.964u_x + 0.244u_{xx} \), with an average coefficient error of only 0.0080. Under high noise conditions (noise level 1.6), the identified equation is \( -0.906u_x + 0.233u_{xx} \), with a higher average parameter error of 0.0550. These results demonstrate that the ANN-PYSR algorithm is capable of reliably extracting the main structure of the equations even in noisy environments, underscoring its robustness and resilience to noise.

Moreover, the running time of the algorithm remained consistent across different noise levels, with an average computation time of approximately $50$ seconds. This further confirms the algorithm’s efficient computational performance and its suitability for PDE recognition tasks in complex, real-world systems.

\begin{table}[htbp]
\centering
\caption{The complexity, loss, and score of the convection-diffusion equation under 60\% noise level}
\label{tag8}
\begin{adjustbox}{max width=\textwidth}
\begin{tabular}{ccccc}
\toprule
Complexity & Loss & Score & Identified PDE \\
\midrule
1  & 0.327263 & 0.000000 & $u_t = -0.0381$ \\
3  & 0.144907 & 0.407335 & $u_t = -0.9558 \cdot u_{x}$ \\
5  & 0.074656 & 0.331602 & $u_t = -1.0462 \cdot u_{x} - 1.0462 \cdot u$ \\
\textbf{7}  & \textbf{0.013918} & \textbf{0.839841} & $\mathbf{u_t = -0.9176 \cdot u_{x} + 0.2433 \cdot u_{xx}}$ \\
9  & 0.013770 & 0.005375 & $u_t = -0.9038 \cdot u_{x} + 0.2430 \cdot u_{xx} - 0.0122$ \\
11 & 0.012884 & 0.033227 & $u_t = 0.0672 \cdot u_{x}^2 - 0.9114 \cdot u_{x} + 0.2456 \cdot u_{xx}$ \\
13 & 0.012236 & 0.025804 & $u_t = 0.0672 \cdot u_{x}^2 - 0.9090 \cdot u_{x} - 0.0259 + 0.2456 \cdot u_{xx}$ \\
\bottomrule
\end{tabular}
\end{adjustbox}
\end{table}

 \begin{table}[h!]
 \centering
 \caption{Identified the convection-diffusion equation, SR Time, and Coefficient Errors for different Noise Levels}
 \label{tag9}
 \begin{tabular}{cccc}
 \toprule
 Noise Level & Identified PDE & Time of SR (s) & Average Coefficient Error \\
 \midrule
 0.2 & \( u_t = -0.9642 u_x + 0.2441 u_{xx} \) & 49.36 & 0.0080 \\
 0.4 & \( u_t = -0.9445 u_x + 0.2374 u_{xx} \) & 51.47 & 0.0140 \\
 0.6 & \( u_t = -0.9175 u_x + 0.2432 u_{xx} \) & 50.15 & 0.0205 \\
 0.8 & \( u_t = -0.8912 u_x + 0.2391 u_{xx} \) & 49.00 & 0.0345 \\
 1.0 & \( u_t = -0.8483 u_x + 0.2244 u_{xx} \) & 50.72 & 0.0640 \\
 1.2 & \( u_t = -0.9192 u_x + 0.2268 u_{xx} \) & 51.45 & 0.0520 \\
 1.4 & \( u_t = -0.9014 u_x + 0.2112 u_{xx} \) & 49.40 & 0.0695 \\
 1.6 & \( u_t = -0.9067 u_x + 0.2336 u_{xx} \) & 50.41 & 0.0550 \\
 \bottomrule
 \end{tabular}
 \end{table}

\subsection{Chaffee-Infante Equation}

The Chaffee-Infante (CI) equation has broad applications across various fields such as environmental science, fluid dynamics, high-energy physics, and electronics. Its mathematical form is given by
\begin{equation}
    u_t = u_{xx} - u + u^3.
\end{equation}
In this study, The training data in this study derived from 301 observation points in the spatial domain \(x \in [0, 3]\) and 200 observation points in the temporal domain \(t \in (0, 0.5)\), resulting in a total of \(301 \times 200 = 60,200\) data points. To evaluate the reliability of the algorithm in PDE identification and the capture of the dynamic structure, we apply noise levels ranging from 0\% to 80\% to the data.

\begin{figure}
    \centering
    \includegraphics[width=0.8\linewidth]{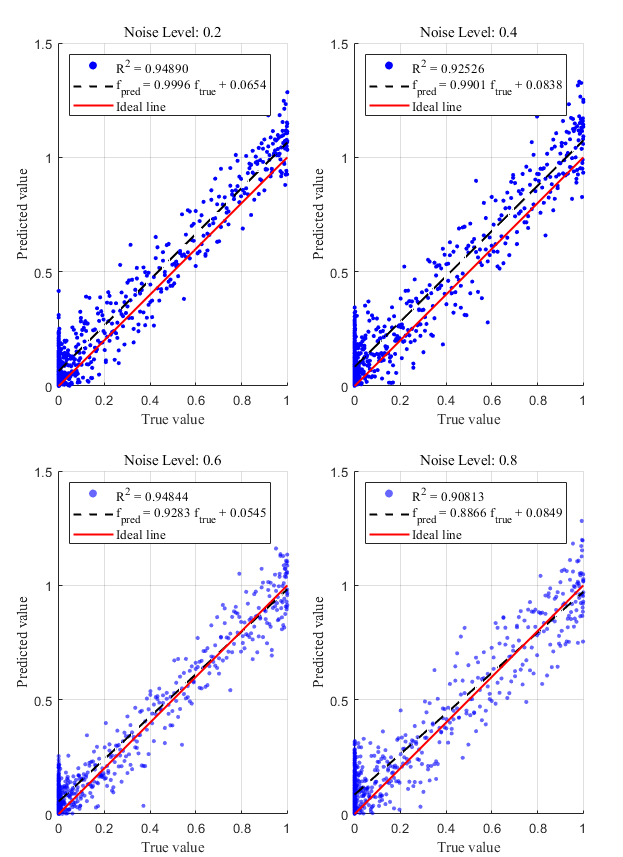}
    \caption{The effect of different noise levels on the prediction accuracy of the CI model. The blue points indicate the scatter distribution of actual data values versus predicted values, the bold line indicates the ideal prediction line, and the black dashed line is the fitted regression line. The corresponding \(R^2\) values and fitted equations are included in the upper right corner of each plot.}
    \label{fig9}
\end{figure}

We examine the identification results of the CI equation at different noise levels (0.2, 0.4, 0.6, and 0.8) using the information presented here. The results, shown in Figure \ref{fig9}, demonstrate that as noise levels increase, the scatter distributions of the predicted and true values become more dispersed. However, the overall linear correlation remains high, despite the increasing deviation between the fitted curves and the ideal line. By calculating the mean squared error and \(R^2\) values, we observe that at lower noise levels (0.2 and 0.4), the model predicts \(R^2\) values of 0.9489 and 0.9252, respectively. As noise levels rise to 0.8, the \(R^2\) values decrease to 0.9081. Since \(R^2\) represents the fraction of variance in the true values explained by the model, where an \(R^2\) of 1 denotes perfect agreement, these findings indicate that the algorithm maintains a high degree of explanatory accuracy even under substantial noise.

At the coefficient level, the prediction results at various noise levels (0.2, 0.4, 0.6, and 0.8) are summarized in Table \ref{tag10}. The average relative errors for these noise levels are 0.179, 0.181, 0.153, and 0.285, respectively. Notably, the coefficient errors for the \(u_{xx}\) and \(u^3\) terms are significantly higher at the highest noise level (0.8), particularly for the \(u^3\) term, which suggests a more pronounced effect of noise on both the linear and nonlinear components of the equation.

\begin{table}[htbp]
\centering
\caption{Chaffee-Infante Equation for Noise Level}
\label{tag10}
\begin{tabular}{cc}
\toprule
Noise Level & Identified PDE \\
\midrule
0.2 & $u_t = 1.108 u_{xx} - 1.312 u + 1.117 u^3$ \\
0.4 & $u_t = 0.908 u_{xx} - 1.212 u + 1.239 u^3$ \\
0.6 & $u_t = 0.893 u_{xx} - 0.788 u + 1.139 u^3$ \\
0.8 & $u_t = 1.393 u_{xx} - 0.824 u + 1.285 u^3$ \\
\bottomrule
\end{tabular}
\end{table}

\subsection{Wave Equation}

Wave equations are a fundamental class of PDEs that describe the propagation of fluctuating phenomena. The two-dimensional wave equation is an extension of the one-dimensional wave equation, used to model wave propagation in two-dimensional media, such as film vibrations, water wave propagation, and the diffusion of acoustic waves in a plane. Its general form is
\begin{equation}
u_{tt} = c^2 (u_{xx} + u_{yy}),
\end{equation}
where \( u(x, y, t) \) represents the displacement at the point \( (x, y) \) in two-dimensional space. The terms \( u_{xx} \) and \( u_{yy} \) describe the propagation of the fluctuation in the respective spatial directions, and \( c \) denotes the wave velocity, which we set to \( c = 1 \) for simplicity. Through the interaction of these two terms, the wave equation exhibits complex fluctuation patterns. As a hyperbolic PDE, the wave equation possesses properties of energy conservation and fluctuation superposition. Identifying the wave equation provides valuable insights into the dynamic properties of the system, shedding light on the mechanisms underlying wave propagation.

\begin{figure}
    \centering
    \includegraphics[width=1\linewidth]{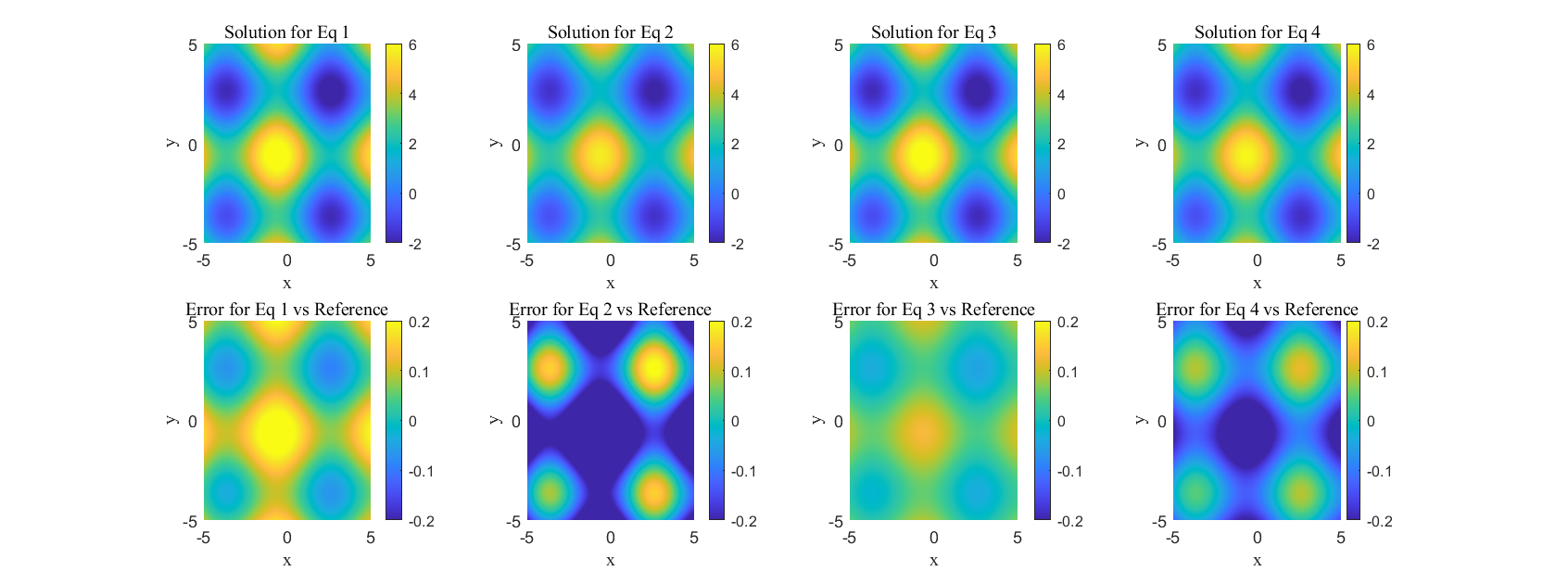}
    \caption{Results of identifying the two-dimensional wave equation \( u_{tt} = u_{xx} + u_{yy} \) using the ANN-PYSR algorithm. The upper part shows the numerical solutions corresponding to the recognized forms of the equation at noise levels of 0.1, 0.2, 0.3, and 0.4, and the lower part shows the error distribution between the numerical solutions of the predicted equation and the numerical solutions of the true equation.}
    \label{fig10}
\end{figure}

The two-dimensional wave equation was identified using the ANN-PYSR algorithm under varying noise conditions. Figure \ref{fig10} illustrates the identification results for the equation \( u_{tt} = u_{xx} + u_{yy} \) at noise levels of 0.1, 0.2, 0.3, and 0.4. The density plot at the top of the figure shows the intensity of the numerical solutions to the discovered equations, while the lower portion displays the error distribution between the true and predicted solutions. Table \ref{tagwave} provides the corresponding predicted equations for each noise level. Overall, the ANN-PYSR algorithm accurately recovers the true equation across all noise levels, with errors primarily concentrated near the boundary regions and at extreme points (peaks and valleys). This indicates that while the accuracy of the algorithm decreases slightly near boundary conditions and at extreme values, the overall error remains within an acceptable range. Additionally, this study also examines the estimation error of the coefficients in the identified equations. When the noise level is 0.1, the coefficient estimation error is minimized to approximately 1.75\%. However, when the noise level increases to 40\%, the estimation error rises to 3.25\%, yet it remains within an acceptable range.

In terms of computational efficiency, the ANN-PYSR algorithm requires an average of 79.21 seconds to identify the two-dimensional wave equation. When compared to traditional numerical methods, the ANN-PYSR algorithm demonstrates superior computational efficiency in handling high-dimensional problems. Despite a slight increase in computation time as the problem complexity increases, the algorithm's efficiency remains favorable. This underscores the practical applicability of the algorithm for real-world systems, capable of identifying complex models within a reasonable time frame.

\begin{table}[htbp]
\centering
\caption{Two-Dimensional Wave Equation for Varying Noise Levels}
\label{tagwave}
\begin{tabular}{cc}
\toprule
Noise Level & Identified PDE \\
\midrule
0.1 & $u_{tt} = 1.0766 u_{xx} + 1.0032 u_{yy}$ \\
0.2 & $u_{tt} = 0.8731 u_{xx} + 0.9522 u_{yy}$ \\
0.3 & $u_{tt} = 0.8354 u_{xx} + 1.2059 u_{yy}$ \\
0.4 & $u_{tt} = 0.7942 u_{xx} + 1.1258 u_{yy}$ \\
\bottomrule
\end{tabular}
\end{table}

\section{Discussion}

This section systematically investigates the performance of the ANN--PYSR algorithm in identifying the Burgers equation (Exp IC) under varying noise levels and sampling densities. Two complementary analyzes, error heatmaps for different noise levels and sample sizes, along with detailed tables of relative errors and identified PDE terms, provide a comprehensive understanding of the robustness of the algorithm.

\begin{figure}
    \centering
    \includegraphics[width=0.6\linewidth]{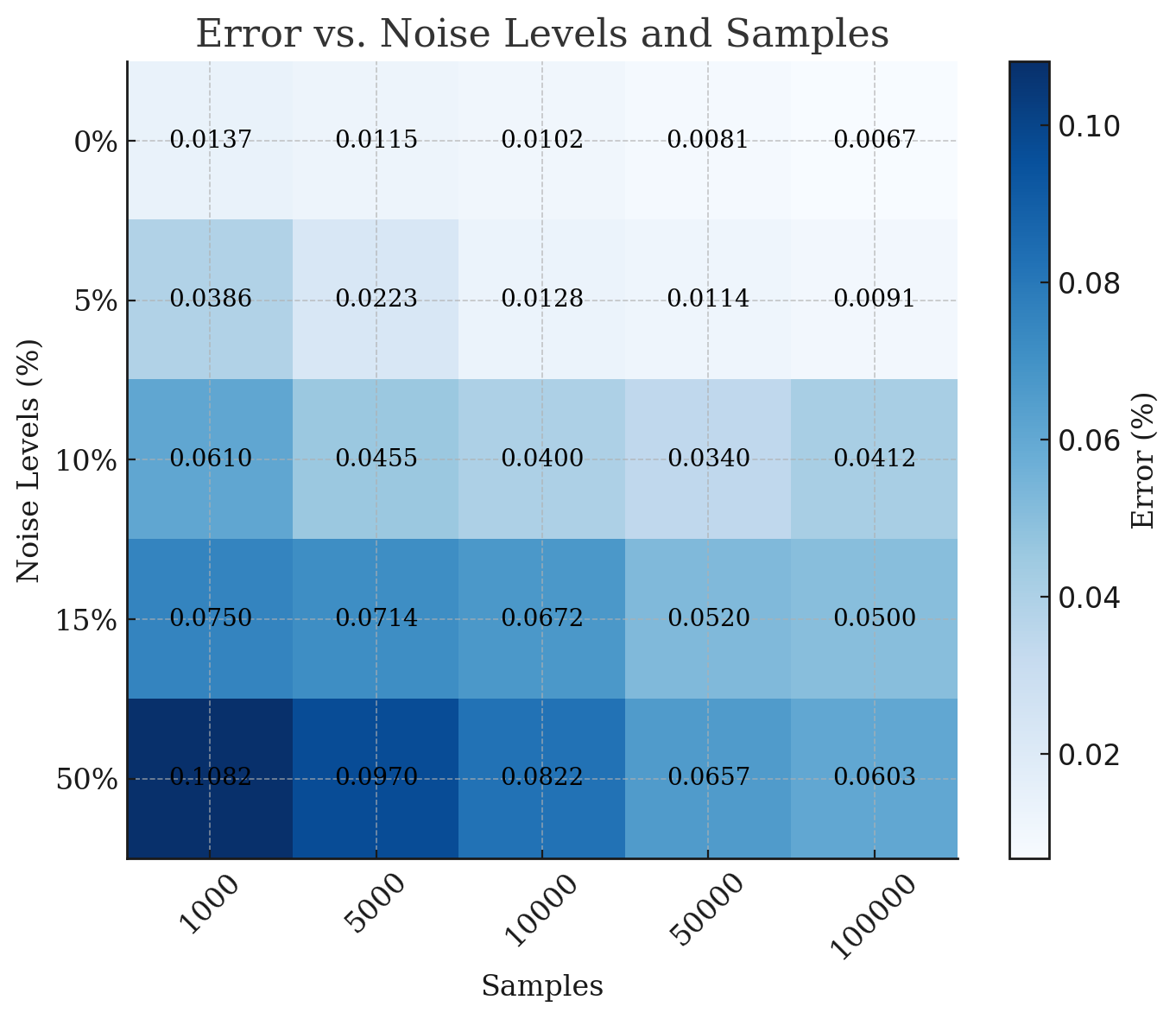}
    \caption{The relationship between noise level, sample size, and the associated error rate. The x-axis represents the number of samples, which ranges from 1,000 to 100,000, while the y-axis represents the noise level, which ranges from 0\% to 50\%. The color scale indicates the percentage error, with lighter shading indicating lower errors and darker shading indicating higher errors.}
    \label{fig11}
\end{figure}

First, the heatmap results (Figure~\ref{fig11}) demonstrate that increasing noise levels (from 0\% to 50\%) leads to higher identification errors, while enlarging the dataset (from 1,000 to 100,000 samples) consistently reduces error. Specifically, at low noise levels (\(\leq\) 10\%), even modest sample sizes (\(\approx\) 10,000) suffice to achieve errors below 5\%. However, as noise levels exceed 30\%, a substantially larger dataset (\(\geq\) 50,000 samples) is necessary to stabilize errors around 10\%. This trade-off underscores the importance of balancing measurement fidelity with data volume when applying data-driven PDE discovery to noisy physical systems. 

Second, the tabulated results for varying noise levels (0.01 to 2.0) and sampling ratios (0.1 to 0.9) reveal how the identified model coefficients evolve, as shown in Table \ref{tag11}. At low noise (\(\leq\) 0.1) and moderate sampling (\(\geq\) 0.5), the recovered equations closely align with the true Burgers dynamics, with relative errors below 2\%. As noise increases, the algorithm compensates by assigning larger weights to the diffusion term \( u_{xx} \), and in extreme cases, it introduces spurious higher-order derivatives. Nevertheless, even with noise = 1.0 and a sampling ratio = 0.3, the core advective-diffusive structure remains discernible, with errors remaining below 15\%.

Taken together, these findings confirm that the ANN–PYSR algorithm can reliably recover key spatio-temporal dynamics from sparse, noisy data, conditions that are typical of real-world sensor networks. In practice, we recommend a sampling density of at least 50\% and noise levels below 20\% to ensure errors remain below 10\%. When higher noise levels are unavoidable, increasing the total number of observations can partially mitigate the degradation in accuracy. Looking ahead, we aim to extend this framework to identify PDEs in more complex, higher-dimensional systems. Furthermore, real-world measurements often exhibit spatial or temporal correlations that are not captured by our current additive white noise model. Future work will explore adaptive sampling strategies, where measurements are concentrated in regions of high solution variability, and incorporate additional physical priors (e.g. conservation laws, symmetries) into the learning process. These enhancements promise to further improve the algorithm’s efficiency and reliability, enabling the robust discovery of more intricate dynamical systems in noisy, data-scarce environments.

\begin{table}[htbp]
\centering
\caption{Relative Errors and Identified Equations for Different Noise Levels and Sampling Ratios}
\label{tag11}
\begin{tabular}{ccc l}
\hline
\textbf{\footnotesize Noise Level} & \textbf{\footnotesize Sampling Ratio} & \textbf{\footnotesize Average Relative Error} & \textbf{\footnotesize \quad \quad \quad \quad Identified Equation} \\ \hline
\small
0.01 & 0.1 & 0.0884 & $u_t = -0.9580 \cdot u \cdot u_{x} + 0.0865 \cdot u_{xx}$ \\
0.05 & 0.1 & 0.0535 & $u_t = -0.9539 \cdot u \cdot u_{x} + 0.0939 \cdot u_{xx}$ \\
0.1 & 0.1 & 0.0556 & $u_t = -0.9568 \cdot u \cdot u_{x} + 0.0931 \cdot u_{xx}$ \\
0.2 & 0.1 & 0.0761 & $u_t = -0.9458 \cdot u \cdot u_{x} + 0.0901 \cdot u_{xx}$ \\
0.1 & 0.3 & 0.0168 & $u_t = -0.9981 \cdot u \cdot u_{x} + 0.0968 \cdot u_{xx}$ \\
0.2 & 0.3 & 0.0111 & $u_t = -0.9962 \cdot u \cdot u_{x} + 0.0981 \cdot u_{xx}$ \\
0.3 & 0.3 & 0.0066 & $u_t = -1.0012 \cdot u \cdot u_{x} + 0.0988 \cdot u_{xx}$ \\
0.4 & 0.3 & 0.0073 & $u_t = -1.0068 \cdot u \cdot u_{x} + 0.0992 \cdot u_{xx}$ \\
0.5 & 0.3 & 0.0163 & $u_t = -1.0168 \cdot u \cdot u_{x} + 0.0984 \cdot u_{xx}$ \\
0.6 & 0.3 & 0.0199 & $u_t = -1.0107 \cdot u \cdot u_{x} + 0.0971 \cdot u_{xx}$ \\
0.7 & 0.4 & 0.0539 & $u_t = -1.0056 \cdot u \cdot u_{x} + 0.0897 \cdot u_{xx}$ \\
0.8 & 0.4 & 0.2169 & $u_t = -1.0192 \cdot u \cdot u_{x} + 0.1414 \cdot u_{xx}$ \\
0.8 & 0.5 & 0.0617 & $u_t = -0.9783 \cdot u \cdot u_{x} + 0.0898 \cdot u_{xx}$ \\
0.9 & 0.5 & 0.0380 & $u_t = -0.9801 \cdot u \cdot u_{x} + 0.0943 \cdot u_{xx}$ \\
1.0 & 0.5 & 0.0738 & $u_t = -0.9762 \cdot u \cdot u_{x} + 0.0872 \cdot u_{xx}$ \\
1.1 & 0.5 & 0.1268 & $u_t = -0.9770 \cdot u \cdot u_{x} + 0.1230 \cdot u_{xx}$ \\
1.1 & 0.7 & 0.0549 & $u_t = -0.9405 \cdot u \cdot u_{x} + 0.1050 \cdot u_{xx}$ \\
1.2 & 0.7 & 0.0512 & $u_t = -0.9387 \cdot u \cdot u_{x} + 0.1041 \cdot u_{xx}$ \\
1.4 & 0.7 & 0.0455 & $u_t = -0.9361 \cdot u \cdot u_{x} + 0.0972 \cdot u_{xx}$ \\
1.5 & 0.8 & 0.0556 & $u_t = -1.0409 \cdot u \cdot u_{x} + 0.0929 \cdot u_{xx}$ \\
1.6 & 0.8 & 0.0633 & $u_t = -1.0574 \cdot u \cdot u_{x} + 0.0930 \cdot u_{xx}$ \\
1.7 & 0.8 & 0.0643 & $u_t = -1.0554 \cdot u \cdot u_{x} + 0.0926 \cdot u_{xx}$ \\
1.8 & 0.9 & 0.0855 & $u_t = -0.9625 \cdot u \cdot u_{x} + 0.1133 \cdot u_{xx}$ \\
1.9 & 0.9 & 0.0887 & $u_t = -0.9749 \cdot u \cdot u_{x} + 0.1152 \cdot u_{xx}$ \\
2.0 & 0.9 & 0.3870 & $u_t = -0.9743 \cdot u \cdot u_{x} + 0.1748 \cdot u_{xx}$ \\ \hline
\end{tabular}
\end{table}

\section{Conclusion}

In this paper, we present an innovative algorithmic framework, ANN-PYSR, which combines symbolic regression and deep learning to automatically identify PDEs from noisy and sparsely sampled observational data. To enhance the accuracy and robustness of the algorithm, we first apply the Savitzky-Golay filter to denoise the data, reducing interference from noise. Next, we employ a residual attention mechanism network to generate high-precision derivative solutions, further improving the model's ability to capture high-frequency features. Finally, the PYSR symbolic regression algorithm is used to identify the governing PDEs of the system. Extensive validation on various classic PDEs, such as the Burgers, Chaffee-Infante, and Korteweg-de Vries equations, shows that ANN-PYSR can accurately capture the primary dynamical structure of the system under extreme conditions, with noise levels up to 200\% and as few as 1000 data points. Compared to traditional methods, ANN-PYSR demonstrates superior precision and computational efficiency when processing high-noise, sparse sampling data. This makes it especially beneficial for sparse sensor networks. Its robustness is further shown in solving complex PDEs. For example, in identifying the Burgers equation, ANN-PYSR achieved an average coefficient error of 0.0015 and a computation time of 47.86 seconds. This significantly outperformed DLGA (error: 0.4000, time: 1973.32 seconds) and R-DLGA (error: 0.0725, time: 2162.71 seconds). In the case of the KdV equation, ANN-PYSR produced an error of 0.0085 with a computation time of 150.29 seconds, surpassing both DLGA and NN-GA. These results emphasize the remarkable accuracy and efficiency of ANN-PYSR, demonstrating its strong potential in noisy and sparse data environments. The goal of ANN-PYSR is to provide more reliable support for modeling in high-noise environments, thereby expanding its applicability in real-world scenarios.

\section*{Acknowledgments}
This work was supported by NSFC Project (12431014) and Project of Scientiﬁc Research Fund of the Hunan Provincial Science and Technology Department (2024ZL5017).


\end{document}